\documentclass[12pt]{article}
\usepackage{graphicx} 

\usepackage{amsmath}
\usepackage[colorlinks=true, allcolors=blue]{hyperref}
\usepackage{subfig}
\usepackage{amssymb}
\usepackage{bm}
\usepackage{float} 
\usepackage{color, xcolor}
     
\usepackage{algorithm}
\usepackage{algorithmic}
\usepackage[margin=1in]{geometry}

\usepackage{makecell}

\providecommand{\keywords}[1]
{
  \small	
  \textbf{\textit{Keywords---}} #1
}

\usepackage{multirow} 
\usepackage{booktabs}
\usepackage{makecell}
\usepackage{changes}

\title{Physics-Informed Learning of Probabilistic Gegenbauer Reconstruction for Transport-Dominated Problems
}

\author{Lei Yan
    \thanks{School of Mathematical Sciences, University of Science and Technology of China, Hefei, Anhui 230026, China.
    E-mail: {\tt lei\_yan@mail.ustc.edu.cn}.}
\and Yan Jiang
    \thanks{School of Mathematical Sciences, University of Science and Technology of China, Hefei, Anhui 230026, China.
    Email: {\tt jiangy@ustc.edu.cn}.
    Research supported by NSFC grant 12271499.}
    }
\date{}

\begin{document}
    \maketitle

    \begin{abstract}

        Transport-dominated problems remain challenging for data-driven methods, which often exhibit severe numerical oscillations in the presence of shocks or steep gradients. These oscillations primarily arise from the use of globally supported basis functions or overly smooth hypothesis spaces that fail to accurately capture localized transport structures.
        \emph{Gegenbauer reconstruction}, originally developed to suppress the Gibbs phenomenon in spectral approximations, has recently been extended to reduced-order modeling and demonstrated promising capability in mitigating such spurious oscillations \cite{yan2026gegenbauer}. Its effectiveness, however, depends critically on the choice of the reconstruction parameters, particularly the weight parameter $\lambda$ and the truncation order $m$. In data-driven applications, the diversity of governing problems, training data, and model architectures makes the design of a unified and theoretically justified parameter-selection strategy highly challenging.

        To address this issue, we propose a physics-informed machine-learning framework that predicts a probability distribution over candidate Gegenbauer parameter pairs, enabling a probabilistically weighted reconstruction that explicitly accounts for parameter uncertainty.
        To balance computational cost and reconstruction accuracy, a two-stage training strategy is adopted, where a general predictor is first pre-trained and subsequently fine-tuned for target problem instances. 
        The proposed framework is evaluated on two representative classes of data-driven solvers, namely reduced-order models and neural operator models, with the POD-Galerkin method and Deep Operator Network (DeepONet) serving as canonical examples, respectively. 
        Numerical experiments on one- and two-dimensional transport-dominated problems demonstrate that the proposed framework successfully learns effective spatially adaptive distributions of Gegenbauer parameters. 
        Compared with conventional reconstruction strategies, the proposed weighted reconstruction reduces numerical errors by up to one to two orders of magnitude while providing a substantially more favorable accuracy--computational cost trade-off than retraining models for individual target problems.

    \keywords{physics-informed machine-learning, pre-training, fine-tuning, Gegenbauer reconstruction, data-driven, transport-dominated, cost-accuracy trade-off}

    \end{abstract}

    \section{Introduction}

This paper introduces a physics-informed machine-learning framework for post-processing oscillatory data-driven solutions of transport-dominated problems.
It is well known that numerical solutions of conservation laws may develop shocks and other complex non-smooth structures, which has motivated decades of research on non-oscillatory numerical methods~\cite{hesthaven2017numerical}. 
Such features pose significant challenges for data-driven approaches, as these methods often rely on globally supported basis functions or excessively smooth hypothesis spaces that are struggling to accurately represent localized transport features~\cite{peherstorfer2022breaking, abbasi2025challenges}.

Various stabilization strategies have been developed to mitigate oscillations in data-driven approximations of transport-dominated problems.
Within the reduced-order modeling (ROM) framework, commonly adopted approaches include modal filtering~\cite{farcas2022filtering} and the introduction of artificial viscosity~\cite{siena2025stabilized}, both of which aim to damp high-frequency components induced by global basis representations and aggressive modal truncation.
Another line of research seeks to modify the ROM formulation itself to better accommodate the intrinsic translational nature of transport phenomena, as in shifted Proper Orthogonal Decomposition~\cite{reiss2018shifted} and shifted Operator Inference~\cite{issan2023predicting}.
In parallel, significant efforts have been devoted to addressing similar challenges within machine-learning-based solvers. 
Representative approaches include domain decomposition techniques~\cite{jagtap2020conservative, jagtap2020extended}, adaptive sampling\cite{wu2023comprehensive, tang2023pinns}, shock-aware loss functions~\cite{roohi2025shock}, and hybrid methods that incorporate physics-based regularization or artificial dissipation into the learning process~\cite{lanthaler2022nonlinear}. 
Despite these advances, developing a unified and computationally efficient correction strategy that reduces spurious oscillations while preserving sharp transport features remains a challenging open problem.

Post-processing strategies provide a non-intrusive way to correct oscillatory data-driven solutions without modifying the underlying reduced-order or learning model.
Among these, Gegenbauer reconstruction has emerged as a promising approach for mitigating unphysical oscillations in transport-dominated problems~\cite{gottlieb1997gibbs, hesthaven2007spectral}. 
Originally developed to alleviate the Gibbs phenomenon in spectral approximations, it has been successfully extended to the reduced-order modeling framework in our previous work \cite{yan2026gegenbauer}, demonstrating the ability to recover sharp solution features while suppressing spurious oscillations near discontinuities. 
Its non-intrusive nature makes it particularly attractive as a post-processing step, allowing existing ROMs or data-driven solvers to benefit from enhanced solution fidelity without modifying the underlying projection or learning framework.
Despite these advantages, the performance of Gegenbauer reconstruction depends sensitively on the reconstruction parameters, particularly the weight parameter $\lambda$ and the truncation order $m$.
While parameter selection in classical spectral methods is supported by established theoretical guidelines~\cite{gelb2004parameter, gelb2005determining}, such guidance is largely absent in data-driven  settings due to the wide variability of problem configurations. 
As a result, the effectiveness of Gegenbauer reconstruction in ROMs and machine-learning-based models can vary substantially, underscoring the need for systematic and automated strategies for parameter prediction.
The goal of this work is therefore to develop an adaptive, data-driven framework that predicts probability distributions over Gegenbauer reconstruction parameters, enabling weighted reconstructions that reduce spurious oscillations while preserving sharp solution features in reduced-order and machine-learning-based models.

In recent years, physics-informed machine learning (PIML) has become an important paradigm for integrating physical knowledge with data-driven modeling~\cite{quarteroni2025combining}. 
By incorporating governing equations, physical constraints, or reference solution information into the learning process, PIML methods can improve data efficiency and promote consistency with known physical laws~\cite{raissi2019physics,karniadakis2021physics}. 
PIML has been applied to a wide range of problems, including fluid dynamics, transport phenomena, and inverse problems, where it has shown improved generalization and stability over purely data-driven approaches in many settings~\cite{kissas2020machine,suzuki2023physics,he2021physics}.  
To improve computational efficiency and training robustness, pre-training and fine-tuning strategies have also been incorporated into physics-informed learning frameworks. 
Originating from transfer learning, this paradigm first learns reusable representations from broad offline data and then adapts them to downstream tasks using a smaller task-specific dataset, thereby balancing computational cost and predictive accuracy~\cite{krizhevsky2012imagenet,tay2021scale,yosinski2014transferable}. 

In this paper, we propose a physics-informed machine-learning framework for probabilistic prediction of Gegenbauer reconstruction parameters from oscillatory data-driven solutions.
The model learns an input-dependent probability distribution over candidate parameter pairs from solution data and reconstruction-induced error information. 
The predicted distribution is then used to form a weighted Gegenbauer reconstruction, which accounts for parameter uncertainty and reduces spurious oscillations while preserving sharp solution features.
To balance reconstruction accuracy and computational cost, we adopt a two-stage training strategy: pre-training on broad oscillatory datasets to learn transferable parameter-prediction patterns, followed by fine-tuning on small target-specific datasets to adapt the predictor to problem-dependent solution features.
Numerical results show that this strategy achieves reconstruction accuracy close to full re-training while requiring substantially lower data-generation and training costs.

The rest of this paper is organized as follows. 
Section~\ref{reviews} reviews the data-driven models considered in this study, including reduced-order models and machine-learning-based solvers. 
Section~\ref{oscillation_discrimination} introduces a quantitative criterion for discriminating oscillation types.
Section~\ref{methodology} presents the proposed physics-informed learning framework for probabilistic Gegenbauer reconstruction in detail. 
In Section~\ref{numerical_results}, numerical experiments on one- and two-dimensional transport-dominated problems are conducted to evaluate the performance of the proposed method, together with an analysis of the cost–accuracy trade-off as well as robustness and generalization.
Finally, Section~\ref{conclusion} summarizes the conclusions and outlines directions for future research.
    \section{Review of data-driven models \label{reviews}}
This section presents the general formulation of the problems considered in this work, followed by a concise overview of the data-driven modeling approaches used in our tests.  We emphasize that data-driven modeling is a vast research area; our goal here is not to provide a comprehensive survey, but rather to summarize the specific reduced-order and machine-learning-based models relevant to our study.

We consider parametric ordinary differential equations (ODEs) posed on a time interval $[t_i,t_f]$,
\begin{equation}
    \frac{d \bm q}{d t}(t;\mu) = \bm f(\bm q,t;\mu), 
    \qquad 
    \bm q(t_i;\mu) = \bm q_{0},
    \label{problem_setup}
\end{equation}
where $\bm q(t;\mu)\in\mathbb{R}^{n}$ denotes the state vector, 
$\mu\in\mathcal{D}\subset\mathbb{R}$ is a parameter, 
$\bm q_0$ is an initial condition, 
and $\bm f:\mathbb{R}^{n}\times [t_i,t_f]\times\mathcal{D}\to\mathbb{R}^{n}$ is a (possibly nonlinear) evolution operator.  
Such ODE systems may arise directly, or from spatial discretizations of partial differential equations (PDEs).  
In this work, we focus on transport-dominated problems whose solutions may contain shocks, leading to significant modeling challenges for both traditional and data-driven approaches.

\subsection{Reduced-order models}

Reduced-order models aim to reduce the computational cost of parametric dynamical systems by approximating the high-dimensional state in a low-dimensional trial space.
Most ROM frameworks employ an \emph{offline-online} strategy: reduced bases and operators are constructed from snapshots offline, and the reduced system is then evaluated efficiently online.
The development of reduced-order models has been rapid, spanning both linear-subspace projection methods and nonlinear-manifold approaches.  
Typical representatives of the former include POD-Galerkin~\cite{hesthaven2016certified} and Operator Inference (OpInf)~\cite{peherstorfer2016data}, while a representative method in the latter category is CAE-LSTM~\cite{maulik2021reduced}.  
In this work, we use the POD-Galerkin reduced-order model, referred to as G-ROM, as the representative ROM setting.

POD~\cite{berkooz1993proper} extracts dominant coherent structures from solution snapshots and, combined with Galerkin projection, provides a classical framework for reduced-order modeling of dynamical systems.
The POD-Galerkin reduced-order model starts from a collection of high-fidelity solution snapshots. 
Given parameter samples $\{\mu_j\}_{j=1}^{n_\mu}\subset \mathcal D$ and time instances $\{t_i\}_{i=1}^{n_t}$, the snapshot matrix is assembled as
\begin{equation}
    \bm{\mathcal Q}
    =
    [\bm q(t_i;\mu_j)]_{i=1,\ldots,n_t;\,j=1,\ldots,n_\mu}
    \in \mathbb R^{n\times n_t n_\mu}.
\end{equation}
The leading $r$ POD modes $\{\bm\varphi_i\}_{i=1}^r$ are obtained from the dominant left singular vectors of $\bm{\mathcal Q}$, equivalently from the eigenvectors of $\bm{\mathcal Q}\bm{\mathcal Q}^T$ associated with the largest eigenvalues. 
The high-dimensional state is then approximated in the reduced space as
\begin{equation}
    \bm q_r(t;\mu)
    =
    \sum_{i=1}^{r} \hat q_i(t;\mu)\bm\varphi_i,
\end{equation}
where $\hat q_i(t;\mu)$ denote the reduced coordinates. 
Applying Galerkin projection to \eqref{problem_setup} yields the reduced dynamical system
\begin{equation}
    \frac{d\hat q_i(t;\mu)}{dt}
    =
    \bigl(\bm f(\bm q_r,t;\mu),\bm\varphi_i\bigr),
    \qquad i=1,\ldots,r,
    \label{eq_GROM}
\end{equation}
with the initial condition $\hat{q}_i(0) = (\bm q_0,\bm \varphi_i)$, where $(\cdot,\cdot)$ denotes the Euclidean inner product.

Since $r\ll n$, only a small number of reduced coordinates need to be evolved online, leading to substantial computational savings compared with the high-fidelity model. 
For nonlinear systems, hyper-reduction techniques such as the Discrete Empirical Interpolation Method (DEIM) \cite{chaturantabut2010nonlinear} can be incorporated to further reduce the online cost. 
However, for convection- or transport-dominated problems, the use of global POD modes and modal truncation may lead to nonphysical oscillations near shocks or steep gradients. 
This motivates stabilization and post-processing strategies, including filtering, spectral-viscosity methods \cite{sirisup2004spectral,bergmann2009enablers,wells2017evolve}, and the Gegenbauer-based reconstruction framework considered in this work.

\subsection{Machine-learning-based models}

Data-driven approximation of operators between function spaces provides an efficient way to construct surrogates that bypass the need for repeatedly solving PDEs during many-query tasks. 
After an offline stage of data generation and training, such surrogates offer rapid online evaluations and have therefore become central tools in scientific machine learning. 
Among the wide range of neural-operator and neural-PDE approaches, physics-informed neural networks (PINNs) \cite{raissi2019physics}, Fourier neural operators (FNOs)\cite{li2020fourier}, and Deep Operator Networks (DeepONets)\cite{lu2019deeponet} are practical and widely-used. 
In the following, we focus on DeepONets as a representative model.

The DeepONet architecture is theoretically motivated by the universal approximation theorem for nonlinear operators \cite{chen1995universal}, which shows that continuous operators between suitable function spaces can be approximated by neural-network-based representations. 
This result naturally leads to the branch-trunk decomposition used in DeepONets: the branch network $\{b_k(u)\}$ encodes the input function into latent coefficients, while the trunk network $\{t_k(x)\}$ learns coordinate-dependent basis functions. 
Accordingly, an operator $\mathcal{G}: u \mapsto \mathcal{G}(u)$ is represented as
\begin{equation}
    \mathcal{G}(u)(x)
    = \sum_{k=1}^{p} b_k(u)\, t_k(x),
\end{equation}
where $p$ denotes the dimension of the latent representation. Both the branch and trunck networks are typically parameterized by feed-forward architectures and are trained jointly on input-output function pairs.

A practically important component of DeepONet is the use of \emph{feature expansion} or \emph{input encoding}. 
Since the input function and output coordinate are represented by finite samples or features, enriched encodings can improve the representation of multiscale and oscillatory structures. 
For periodic or wave-like problems, Fourier-type encodings are commonly used, for example
\begin{equation}
    \label{eq:deeponet_periodic_encoding}
    x \;\mapsto\; 
    \bigl[
    x,\;
    \sin(\pi x),\; \cos(\pi x),\;
    \sin(2\pi x),\; \cos(2\pi x),\;
    \dots,\;
    \sin(K\pi x),\; \cos(K\pi x)
    \bigr].
\end{equation}
which is particularly useful for transport- and advection-dominated solution operators.

Several DeepONet variants have been proposed to further improve efficiency and generalization: POD-DeepONet\cite{lu2022comprehensive} replaces the trunk net with the modes calculated on training data by POD, aligning the learned representation with the dominant modes of the PDE solution manifold and reducing the number of trainable parameters; 
DON-LSTM \cite{michalowska2025multi} extends the DeepONet with a long short-term memory network (LSTM) to leverage multi-resolution data; 
Physics-informed DeepONet\cite{goswami2022physics} embeds PDE residuals directly into the loss function. 
Although differing in architectural details, these variants adhere to the same operator-separation principle that underlies the original DeepONet framework.
    
\section{Oscillation Discrimination}
\label{oscillation_discrimination}

The Gegenbauer reconstruction is performed independently over each analyticity interval bounded by two adjacent discontinuities. Therefore, accurate discontinuity detection is an essential prerequisite.
In one-dimensional problems, Gibbs-type oscillations, when present, are typically localized near genuine discontinuities and can therefore serve as useful indicators of interval boundaries. However, this strategy is not directly applicable in two dimensions. 
Since the proposed Gegenbauer reconstruction is implemented in a dimension-by-dimension manner, the observed oscillations along a one-dimensional cut may either be associated with a discontinuity in the current direction or arise from the influence of a nearby discontinuity in the transverse direction, even though the solution is smooth in the current direction. 
It is therefore necessary to distinguish \emph{discontinuity-induced oscillations} (DIOs) from \emph{spurious oscillations in smooth regions} (SOSs).

In the present work, candidate discontinuity locations are first identified from the one-dimensional profile using a first-derivative Sobel-type gradient operator~\cite{sobel19683x3}. Let $\mathcal S^{(1)}$ denotes the corresponding discrete first-derivative operator, and for the grid values $u_i=u(x_i)$, the local gradient indicator is defined as
\begin{equation}
    \label{eq:sobel_indicator}
    G_i = \left| (\mathcal S^{(1)}u)_i \right|.
\end{equation}
Large values of $G_i$ indicate sharp local variations and are therefore used to generate candidate discontinuity locations.  Grid points satisfying $G_i>\tau_{\mathrm S}$ are marked as candidates, where $\tau_{\mathrm S}$ is a prescribed threshold. 
Neighboring candidates are clustered, and the point with the largest $G_i$ in each cluster is retained as the representative discontinuity location. 
After sorting the detected locations $\{s_j\}_{j=1}^{K}$, the whole computational domain $\Omega=[x_{\min},x_{\max}]$ is partitioned into candidate analyticity intervals, with $[a,b]$ denoting one such interval for the local Gegenbauer reconstruction.
The DIO/SOS classification then distinguishes boundary-localized oscillations associated with genuine discontinuities from spurious high-frequency oscillations occurring within otherwise smooth intervals.
To this end, we introduce two complementary criteria based on the spatial distribution of high-frequency energy.

Given the solution $u(x)$ on the analyticity interval $[a,b]$, high-frequency components are isolated using a bandpass filter~\cite{oppenheim1999discrete}
\begin{equation}
    u_h(x) = \mathcal{F}^{-1}\!\left[ H(f)\, \mathcal{F}[u(x)] \right],
\end{equation}
where $\mathcal{F}$ denotes the Fourier transform and $H(f)$ retains frequencies within a prescribed range.
Let $\{x_i\}$ denotes the uniform grid points in the analyticity interval $[a,b]$, the dominant oscillation location is characterized by the energy-weighted coordinate
\begin{equation}
    x_c = \frac{\sum_i x_i e(x_i)}{\sum_i e(x_i)}, 
    \qquad
    x_{\mathrm{norm}} = \frac{x_c - a}{b-a},
\end{equation}
where $e(x) = u_h(x)^2$ denotes high-frequency energy. 
It is observed that the location of $x_c$ denotes the energy centroid of the high-frequency component of $e_h$, while the normalized quantity $x_{\mathrm{norm}}$ characterizes the spatial distribution of the high-frequency oscillatory energy within the analyticity interval. 
When the high-frequency energy is concentrated near the boundaries of the analyticity interval, that is, when $x_{\mathrm{norm}}$ lies close to 0 or 1, the oscillations are classified as DIOs. Conversely, when the high-frequency energy is concentrated in the interior of the interval, that is, when $x_{\mathrm{norm}}$ lies away from the 0 and 1, the oscillations are classified as SOSs.

To make this criterion quantitative, we partition the analyticity interval $[a,b]$ into three regions: 
the left-edge region $\mathcal{I}_L = [a,a+\frac{1}{4}(b-a))$, the interior region $\mathcal{I}_M = [ a+\frac{1}{4}(b-a), b-\frac{1}{4}(b-a)]$, and the right-edge region $\mathcal{I}_R=(b-\frac{1}{4}(b-a), b]$.
The average high-frequency energies over these regions are
\begin{equation}
    E_\ell = \frac{1}{|\mathcal{I}_\ell|} \sum_{x_i\in I_\ell} e_h(x_i),
    \qquad
    \ell \in \{L,M,R\}.
    \label{eq:regional_energy}
\end{equation}
The edge-to-interior energy ratio is then defined as
\begin{equation}
    S = \frac{\max(E_L,E_R)}{E_M+\delta},
    \label{eq:edge_energy_ratio}
\end{equation}
where $\delta=10^{-14}$ is a small regularization constant to avoid division by zero.
A large value of $S$ confirms that the high-frequency energy is mainly concentrated near the interval boundaries, whereas a small one suggests the oscillation is distributed in the interior.

The indicators $x_{\rm norm}$ and $S$ provide complementary information: the former characterizes the location of the high-frequency energy, whereas the latter quantifies its edge-to-interior distribution.
Finally, combining these two indicators, the oscillation type is classified according to
\begin{equation}
    \text{oscillation type} =
    \begin{cases}
    \text{DIOs}, & 
    \begin{aligned}[t]
    & S > S_{\mathrm{thresh}} \;\text{or}\;
    x_{\mathrm{norm}} < \mathrm{edge\_frac} \\
    & \text{or}\;
    x_{\mathrm{norm}} > 1-\mathrm{edge\_frac},
    \end{aligned}
    \\[1mm]
    \text{SOSs}, & \text{otherwise},
    \end{cases}
\end{equation}
where $S_\mathrm{thresh}$ is chosen around $2\sim 3$, and $\mathrm{edge\_frac}$ around $0.2\sim 0.25$ in the present experiments. 
This criterion provides a problem-dependent working definition for separating boundary-localized, discontinuity-induced oscillations from interior spurious oscillations, and is used to support the subsequent dataset construction and training procedure.

    \section{Methodology \label{methodology}}

The objective of this work is to develop an automatic parameter-prediction model for Gegenbauer reconstruction to mitigate oscillations in transport-dominated problems. 
This section details the proposed physics-informed probabilistic model, including motivation, model design and numetrical realization.

\subsection{Gegenbauer reconstruction}
\label{Gegenbauer_recon}

Oscillations in spectral approximations of discontinuous solutions, commonly known as the Gibbs phenomenon, are well documented. 
The \emph{Gegenbauer reconstruction} provides a classical post-processing strategy for alleviating such oscillations by re-expanding a truncated global approximation in a Gegenbauer basis on analyticity subintervals~\cite{gottlieb1997gibbs}.

Specifically, starting from the truncated spectral expansion
\begin{equation}
    f_N(x) = \sum_{k=-N}^{N} (f,\Psi_k)\Psi_k(x),
\end{equation}
where $\{\Psi_k\}$ denotes the original spectral basis such as Fourier modes or orthogonal polynomials, the reconstruction is performed locally on each interval where the solution is analytic.

Let $\{C_k^\lambda (\xi)\}_{k\geq 0}$ denote the Gegenbauer polynomials orthogonal on $[-1,1]$ with respect to the weight $(1-\xi^2)^{\lambda-1/2}$. On an analyticity interval $[a,b]$, mapped to $\xi\in[-1,1]$, the Gegenbauer reconstruction takes the form of
\begin{equation}
    f_N^{\lambda,m}(x) = \sum_{k=0}^{m} \hat{g}_\epsilon^\lambda(k) C_k^\lambda(\xi(x)), \qquad 
    \hat{g}_\epsilon^\lambda(k) = \frac{1}{h_k^\lambda}\langle f_N(x(\xi)), C_k^\lambda \rangle_\lambda, \quad k=0,\dots,m,
\end{equation}
where the weighted inner product and the normalization constant $h_k^\lambda$ are defined as
\begin{equation}
    \langle f,g\rangle_{\lambda} = \int_{-1}^{1} (1-\xi^2)^{\lambda-\frac12} f(\xi)\, g(\xi)\, d\xi, \quad 
    h_k^\lambda = \langle C^\lambda_k,C^\lambda_k \rangle_{\lambda}. 
\end{equation}

The convergence and stability of Gegenbauer reconstruction critically depend on the choice of reconstruction parameters $(m, \lambda)$. 
In practice, an inappropriate parameter selection may result in excessive smoothing, loss of sharp features, or amplification of spurious oscillations.
For classical spectral approximations, theoretical results provide guidance for choosing these parameters in relation to the truncation level $N$, leading to spectral accuracy under suitable assumptions~\cite{gottlieb1997gibbs}.

Recent efforts have extended Gegenbauer post-processing to data-driven
models, whose global basis functions often introduce oscillatory artifacts when approximating transport-dominated solutions. 
The methods typically apply the reconstruction on detected analyticity intervals.
However, as discussed in~\cite{yan2026gegenbauer}, unlike classical spectral approximations, data-driven models do not provide a fixed truncation structure from which $(m, \lambda)$ can be selected systematically. 
Instead, the effective approximation space is implicitly shaped by the underlying problem, the training data, the model architecture, and the learning process. Consequently, the optimal reconstruction parameters are intrinsically problem-dependent and may vary across different spatial locations.

These observations motivate the present work, which aims to develop an automated, data-driven parameter-prediction model for Gegenbauer reconstruction. 
In particular, instead of selecting a single parameter pair heuristically, the proposed framework predicts a probability distribution over a set of candidate parameter pairs $(m,\lambda)$.
This probabilistic formulation addresses the instability of hard parameter selection in data-driven reconstruction: for a given oscillatory input profile, several parameter pairs may yield comparable reconstruction errors, while the empirically best pair can vary with the local shock location, oscillation amplitude, model architecture, and training distribution. 
By learning a soft distribution over candidate parameters, the final weighted reconstruction can average over multiple effective candidates, thereby reducing the risk of selecting a suboptimal pair and improving robustness across spatial locations and problem instances.

\subsection{Entropy-Regularized Gegenbauer Parameter Selection}
\label{prob_parameter_prediction}

A remaining question here is how to select an appropriate parameter pair for a given oscillatory data-driven solution. A deterministic minimum-error selection retains only one candidate and may be unstable when several parameter pairs yield comparable reconstruction errors. We therefore formulate the parameter-selection problem as an entropy-regularized probabilistic relaxation over the candidate parameter set.

Let $\mathcal M=\{0,\ldots,m_{\max}\}$ and $\Lambda=\{1,\ldots,\lambda_{\max}\}$ denote the candidate truncation orders and Gegenbauer weight parameters, respectively, and 
Let the Gegenbauer parameter manifold be
\begin{equation*}
    \Gamma
    =
    \mathcal M\times\Lambda
    =
    \{(m_i,\lambda_j):
      m_i \in \mathcal{M};\ 
      \lambda_j\in \Lambda\}.
\end{equation*}

For an oscillatory input profile $v$, denote the Gegenbauer reconstruction associated with the candidate pair $(m_i,\lambda_j)$ as
\begin{equation}
    v_{i,j}(x)
    :=
    \mathcal G_{m_i,\lambda_j}[v](x)
    =
    \sum_{k=0}^{m_i}
    \widehat v_k^{\lambda_j}
    C_k^{\lambda_j}(\xi(x)),
    \label{eq:gegenbauer_recon}
\end{equation}
where $\widehat v_k^{\lambda_j}$ is the Gegenbauer coefficients and $\xi=\xi(x)$ is the coordinate transform.
Given the corresponding reference solution $u_{\mathrm{ref}}$, let reconstruction error associated with $(m_i,\lambda_j)$ be
\begin{equation}
    E_{i,j}(v,u_{\mathrm{ref}})
    =
    \ell\bigl(v_{i,j},u_{\mathrm{ref}}\bigr),
    \label{eq:parameter_reconstruction_error}
\end{equation}
where $\ell$ denotes the prescribed reconstruction-error functional,  possibly evaluated on a masked grid set when small discontinuity-location mismatches are excluded.

The conventional parameter-selection problem seeks a single minimum-error candidate
\begin{equation}
    (i^\star,j^\star)
    \in
    \operatorname*{argmin}_{(m_i, \lambda_j)\in \Gamma }
    E_{i,j}(v,u_{\mathrm{ref}}),
    \label{eq:deterministic_parameter_selection}
\end{equation}
with the minimum-error index set denoted as $\mathcal I_\star = \left\{(i,j): E_{i,j}(v,u_{\mathrm{ref}})=E_{\min} \right\}.$
Define the probability simplex over $\Gamma$ as
\begin{equation}
    \Delta_\Gamma
    =
    \left\{
        P=(P_{i,j})\in\mathbb R^{n_m\times n_\lambda}:
        P_{i,j}\geq 0,\quad
        \sum_{i=1}^{n_m}
        \sum_{j=1}^{n_\lambda}
        P_{i,j}=1
    \right\}
\end{equation}
and the deterministic problem can be then written in the relaxed form
\begin{equation}
    \min_{P\in\Delta_\Gamma}
    \sum_{i=1}^{n_m}
    \sum_{j=1}^{n_\lambda}
    P_{i,j}E_{i,j}(v,u_{\mathrm{ref}}).
    \label{eq:unregularized_probability_relaxation}
\end{equation}
The probability-simplex formulation preserves both the optimal value and the minimizers of the deterministic selection problem. Indeed,
\begin{equation}
    \min_{P\in\Delta_\Gamma}
    \sum_{i=1}^{n_m}\sum_{j=1}^{n_\lambda}
    P_{i,j}E_{i,j}
    =
    \min_{(m_i, \lambda_j)\in \Gamma}
    E_{i,j}
    =
    E_{\min}.
\end{equation}
Moreover, equality is attained if and only if
\begin{equation*}
    P_{i,j}=0
    \qquad
    \text{for all }(i,j)\notin\mathcal I_\star.
\end{equation*}
Thus, if the deterministic minimizer is unique, the unique optimizer over $\Delta_\Gamma$ is the Dirac distribution concentrated on the same minimum-error parameter pair.

To further obtain a smooth and non-degenerate selection rule, we introduce a Shannon entropy regularization~\cite{shannon1948mathematical} over $\Delta_\Gamma$, which is defined as
\begin{equation}
    P_\gamma^\star(v,u_{\mathrm{ref}})
    =
    \operatorname*{argmin}_{P\in\Delta_\Gamma}
    \left\{
        \sum_{i=1}^{n_m}
        \sum_{j=1}^{n_\lambda}
        P_{i,j}E_{i,j}(v,u_{\mathrm{ref}})
        +
        \frac{1}{\gamma}
        \sum_{i=1}^{n_m}
        \sum_{j=1}^{n_\lambda}
        P_{i,j}\log P_{i,j}
    \right\},
    \label{eq:entropy_regularized_selection}
\end{equation}
with the first term favoring parameter pairs with small reconstruction errors, and the entropy term preventing the distribution from prematurely collapsing onto a single candidate.
The optimization problem~\eqref{eq:entropy_regularized_selection} is strictly convex and therefore admits a unique minimizer. Introducing a Lagrange multiplier $\nu$ for the normalization constraint and applying the first-order optimality condition yield
\begin{equation*}
    E_{i,j} + \frac{1}{\gamma} \bigl(1+\log P_{i,j}\bigr) + \nu = 0. 
\end{equation*}
Solving for $P_{i,j}$ and enforcing the normalization constraint give a Gibbs distribution over all candidate parameter pairs
\begin{equation}
    P_{\gamma,i,j}^\star(v,u_{\mathrm{ref}})
    =
    \frac{
        \exp\bigl(
            -\gamma E_{i,j}(v,u_{\mathrm{ref}})
        \bigr)
    }{
        \displaystyle
        \sum_{i'=1}^{n_m}
        \sum_{j'=1}^{n_\lambda}
        \exp\bigl(
            -\gamma E_{i',j'}(v,u_{\mathrm{ref}})
        \bigr)
    },
    \label{eq:gibbs_parameter_distribution}
\end{equation}
where $\gamma>0$ controls the concentration of the probability distribution around the low-error parameter pairs.
This formulation provides a variational justification for converting the reconstruction-error matrix into a probability distribution.
Moreover, it establishes a continuous transition from deterministic to probabilistic parameter selection. 
When the minimum-error pair is unique, $P_\gamma^\star$ converges to the corresponding one-hot distribution as $\gamma\to\infty$, whereas finite values of $\gamma$ retain information from multiple low-error candidates.
The resulting distribution is subsequently used as the target probability law for the learned parameter predictor.

The entropy-regularized formulation introduced above assigns a probability distribution to the candidate Gegenbauer parameter pairs, and we use this distribution to define a weighted reconstruction in the following.
For any $P=(P_{i,j})\in\Delta_\Gamma$, let convex combination of the candidate reconstructions be
\begin{equation}
    \mathcal R(P;v)(x)
    =
    \sum_{i=1}^{n_m}
    \sum_{j=1}^{n_\lambda}
    P_{i,j}\,v_{i,j}(x),
    \label{eq:weighted_reconstruction_operator}
\end{equation}
where $v_{i,j}$ denotes the Gegenbauer reconstruction of the input profile $v$ associated with the parameter pair $(m_i,\lambda_j)$, as defined in Subsection~\ref{Gegenbauer_recon}. 
Indeed, using the entropy-regularized distribution $P_\gamma^\star(v,u_{\mathrm{ref}})$, we obtain the corresponding relaxed reconstruction 
\begin{equation}
    v_\gamma^\star(x)
    =
    \mathcal R
    \bigl(P_\gamma^\star(v,u_{\mathrm{ref}});v\bigr)(x)
    =
    \sum_{i=1}^{n_m}
    \sum_{j=1}^{n_\lambda}
    P_{\gamma,i,j}^\star(v,u_{\mathrm{ref}})
    v_{i,j}(x),
    \label{eq:entropy_relaxed_reconstruction}
\end{equation}
which incorporates contributions from multiple low-error candidate pairs for finite $\gamma$.
However, the distribution $P_\gamma^\star(v,u_{\mathrm{ref}})$ depends on the reference solution and is therefore unavailable during online post-processing. We consequently introduce an input-dependent parametric predictor
\begin{equation}
    P_\theta(v)
    =
    \left\{
        p_\theta(m_i,\lambda_j\mid v)
    \right\}_{
        i=1,\ldots,n_m;\,
        j=1,\ldots,n_\lambda
    }
    \in\Delta_\Gamma,
    \label{eq:predicted_parameter_distribution}
\end{equation}
where $\theta$ denotes the trainable model parameters in neural networks. The predictor is trained to approximate the error-induced selection law $P_\gamma^\star(v,u_{\mathrm{ref}})$ from the oscillatory input profile $v$ alone.
The final learned reconstruction model is then defined as
\begin{equation}
    \widetilde v_\theta(x)
    =
    \mathcal R(P_\theta(v);v)(x)
    =
    \sum_{i=1}^{n_m}
    \sum_{j=1}^{n_\lambda}
    p_\theta(m_i,\lambda_j\mid v)\,
    v_{i,j}(x).
    \label{eq:learned_weighted_reconstruction}
\end{equation}
We emphasize that the predicted distribution $P_\theta(v)$ is not a Bayesian posterior, but an error-induced soft representation of parameter-selection uncertainty. 
Because it directly determines the weighted reconstruction, this uncertainty is an active part of the post-processing procedure rather than an auxiliary output.
In contrast to deterministic selection of a single parameter pair, this probabilistic weighted reconstruction formulation retains information provides a more flexible, robust and input-dependent post-processing mechanism for oscillatory data-driven solutions.

\subsection{\texorpdfstring{Numerical Implementation}{Numerical Implementation}}
\label{training_pipeline}

This subsection outlines the overall training pipeline for the proposed probabilistic Gegenbauer parameter-prediction model. 
The goal of the training process is to learn an input-dependent probability distribution over candidate parameter pairs from oscillatory data-driven solution profiles. 
The workflow consists of three main components: constructing datasets that capture both DIOs and SSOs; training a probabilistic neural model with a reconstruction-informed loss to learn parameter distributions; and a two-stage pre-training and fine-tuning strategy to adapt the model to specific problems, enhancing reconstruction accuracy and efficiency.

\subsubsection{Dataset Generation and Soft-Target Construction\label{dataset}}

This part describes the construction of training datasets and the associated soft probability targets. 
The pre-training dataset is designed to be large and diverse, and therefore emphasizes low-cost sample generation. 
In contrast, the fine-tuning dataset is much smaller and problem-specific, and therefore emphasizes the accuracy of both the data-driven input profiles and the reference solutions. 
After the datasets are generated, the same soft-target construction is applied to all samples by converting the reconstruction-error matrix over candidate Gegenbauer parameters into a probability distribution.

\paragraph{Dataset for pre-training.}
The pre-training dataset is constructed to capture representative oscillatory patterns that commonly appear in data-driven approximations of transport-dominated problems.  
We construct separate pre-training datasets
\begin{equation*}
    \mathcal D_{\mathrm{pre}}^{r}
    =
    \left\{
        \bigl(v^{(s)},u_{\mathrm{ref}}^{(s)}\bigr)
    \right\}_{s=1}^{N_{\mathrm{pre}}^{r}},
    \qquad
    r\in\{\mathrm{DIOs},\mathrm{SSOs}\}
\end{equation*}
for discontinuity-induced oscillations and spurious oscillations in smooth regions, respectively.

For discontinuity-induced oscillations, training samples are generated based on a representative one-dimensional transport problem (see Appendix~\ref{appendix_dataset_true} for detailed setup and parameters). 
To avoid repeatedly solving the full-order or reduced-order dynamical system, shifted smooth profiles are projected onto a precomputed reduced space. The projected profiles retain the localized oscillatory structures produced by truncated global approximations near genuine discontinuities, while requiring substantially lower generation cost than direct simulation.

For spurious oscillations, samples are generated by perturbing a smooth baseline with localized oscillatory and stochastic components,
\begin{equation}
    v_{\mathrm{spu}}(x)
    =
    v_{\mathrm{base}}(x)
    +
    \mathcal O(x)
    +
    \mathcal N(x),
    \label{eq:spurious_training_profile}
\end{equation}
where $v_{\mathrm{base}}$ is a constant or low-degree polynomial, $\mathcal O$ is a localized oscillatory perturbation, and $\mathcal N$ is a localized smooth stochastic perturbation (see Appendix~\ref{appendix_dataset_spurious} for detailed setup and parameters). 
This construction produces controlled nonphysical oscillations without introducing genuine discontinuities.

For both classes, Gegenbauer reconstructions are evaluated over all $(m_i,\lambda_j)\in\Gamma$, and only samples whose best reconstruction error is within an acceptable range are retained.

\paragraph{Dataset for fine-tuning.}

The fine-tuning dataset is generated directly from the target problem, which is intended to adapt the pre-trained parameter predictor to the geometry, amplitude, and error characteristics.
Since only a small number of samples are used for fine-tuning, the dataset is constructed using accurate target-problem information.
Specifically, for a prescribed target PDE and data-driven solver, the fine-tuning dataset is denoted by
\begin{equation*}
    \mathcal D_{\mathrm{tar}}^{r}
    =
    \left\{
        \bigl(v^{(s)},u_{\mathrm{ref}}^{(s)}\bigr)
    \right\}_{s=1}^{N_{\mathrm{tar}}^{r}},
    \qquad
    r\in\{\mathrm{DIOs},\mathrm{SSOs}\},
\end{equation*}
where $v^{(s)}$ is produced by the data-driven model to be post-processed, such as a G-ROM or neural-operator approximation, whereas $u_{\mathrm{ref}}^{(s)}$ is obtained from the corresponding exact or high-fidelity solution. 
When both oscillation types occur, the class index $r$ is determined by the oscillation-discrimination criterion in Section~\ref{oscillation_discrimination}.

\paragraph{Soft-target construction from the error matrix.}
After the input profile and reference solution are obtained, the same soft-target construction is applied to both pre-training and fine-tuning samples. 
For each sample $\bigl(v^{(s)},u_{\mathrm{ref}}^{(s)}\bigr)$ and each $(m_i,\lambda_j)\in\Gamma$, let $v_{i,j}^{(s)} = \mathcal G_{m_i,\lambda_j} \bigl[v^{(s)}\bigr]$ be the Gegenbauer reconstruction associated with parameter pair $(m_i,\lambda_j)\in\Gamma$ and $E_{i,j}^{(s)} =\ell\bigl( v_{i,j}^{(s)}, u_{\mathrm{ref}}^{(s)} \bigr)$ be the reconstruction error, the corresponding reconstruction-error matrix is therefore
\begin{equation*}
    E^{(s)}
    =
    \bigl(E_{i,j}^{(s)}\bigr)
    \in
    \mathbb R^{n_m\times n_\lambda}.
\end{equation*}

The soft target $P_{\gamma}^{\star,(s)}$ associated with the $s$-th sample is the entropy-regularized distribution derived in Subsection~\ref{prob_parameter_prediction}
\begin{equation}
    P_{\gamma,i,j}^{\star,(s)}
    =
    \frac{
        \exp\bigl(-\gamma E_{i,j}^{(s)}\bigr)
    }{
        \displaystyle
        \sum_{i'=1}^{n_m}
        \sum_{j'=1}^{n_\lambda}
        \exp\bigl(-\gamma E_{i',j'}^{(s)}\bigr)
    },
    \qquad
    i=1,\ldots,n_m,\quad
    j=1,\ldots,n_\lambda.
    \label{eq:soft_target}
\end{equation}
Thus, the complete supervised dataset may be written as
\begin{equation*}
    \widehat{\mathcal D}
    =
    \left\{
        \bigl(
            v^{(s)},
            u_{\mathrm{ref}}^{(s)},
            P_{\gamma}^{\star,(s)}
        \bigr)
    \right\}_{s=1}^{N},
\end{equation*}
with lower-error parameter pairs receiving larger probability weights. 
For finite $\gamma$, the target retains information from multiple near-optimal candidates; when the minimum-error pair is unique, it approaches the corresponding one-hot distribution as $\gamma\to\infty$. 
Consequently, $P_{\gamma}^{\star,(s)}$ provides a non-degenerate supervisory target when several parameter pairs yield comparable reconstruction errors.

\subsubsection{Probabilistic Model and Training Objective \label{prob_model_and_loss}}

The goal of the probabilistic model $P_\theta(v)\in \Delta_\Gamma$ is to approximate the error-induced selection mapp redict an input-dependent distribution over the candidate Gegenbauer parameter pairs in $\Gamma$ from the oscillatory input profile $v$
\begin{equation*}
    v
    \longmapsto
    P_\gamma^\star(v,u_{\mathrm{ref}}),
\end{equation*}
with the soft target distribution $P_{\gamma}^{\star,(s)}$ providing the supervisory signal.

\begin{figure}[!ht]
    \centering
    \includegraphics[width=0.95\linewidth]{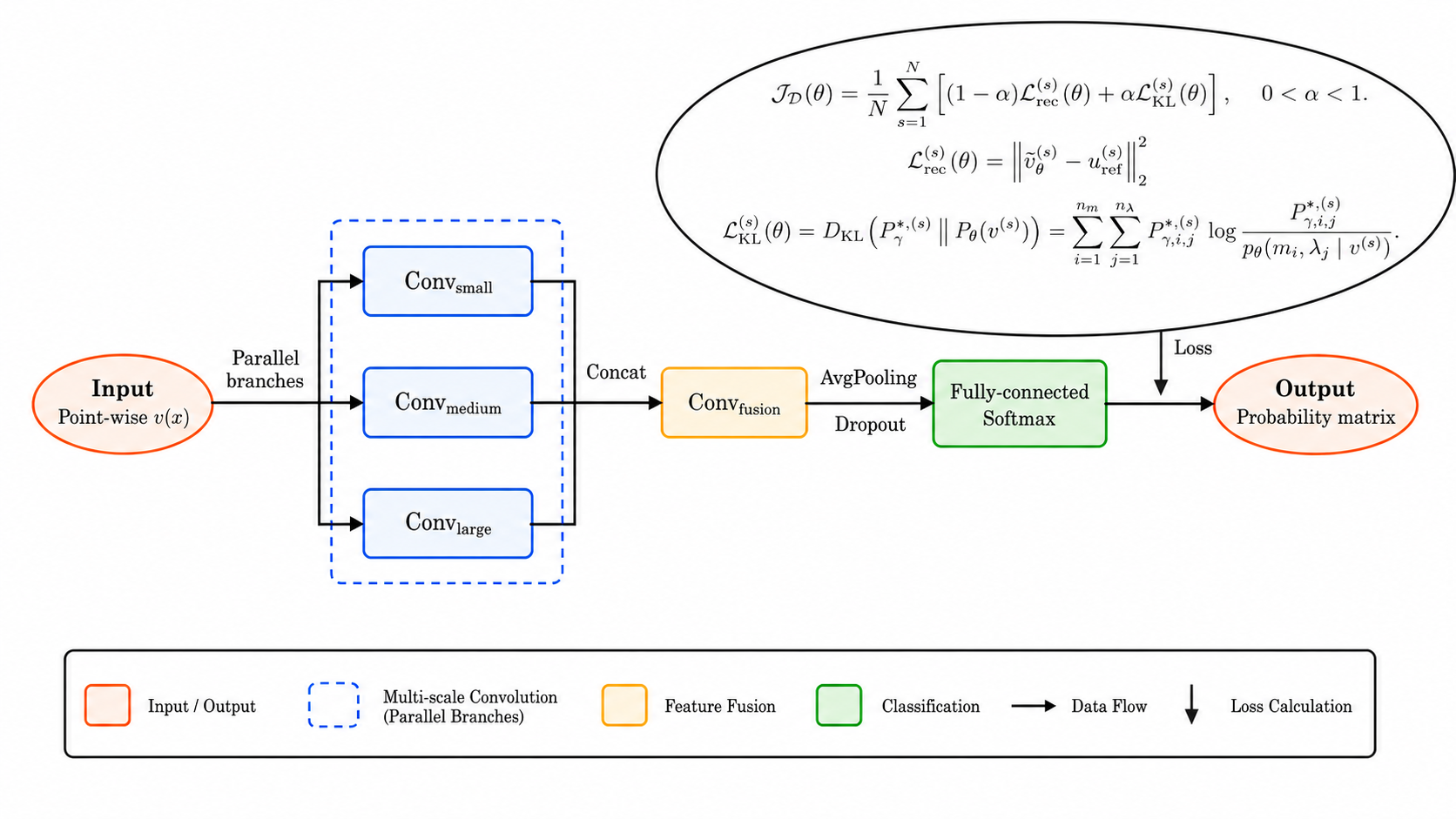}
    \caption{Multi-scale convolutional network architecture for oscillation-aware Gegenbauer parameter inference.
    The input signal is processed by parallel convolutional branches with different receptive-field sizes to extract complementary oscillatory features. A fusion convolution integrates the multi-scale representations into a unified latent feature, which is mapped through a fully-connected layer and a softmax function to produce the predicted parameter distribution. }
    \label{fig_network_architecture}
\end{figure}

\paragraph{Convolutional architecture for oscillation-aware feature extraction.}
The parameter predictor is implemented by a multi-scale convolutional neural network designed to extract oscillation-aware local features. 
Parallel convolutional branches with different kernel sizes provide multiple receptive fields, allowing the network to capture oscillatory patterns with different characteristic spatial scales~\cite{szegedy2015going, luo2016understanding}.
The resulting features are concatenated along the channel dimension, compressed by a fusion convolution, globally averaged over the spatial dimension, and then mapped by a fully connected layer to a $|\Gamma|$-dimensional logits $\{z_{i,j}\}$.
A softmax normalization then gives the predicted probability distribution
\begin{equation}
    p_{\theta}(m_i,\lambda_j \mid v)
    =
    \frac{\exp(z_{i,j})}
    {
    \sum_{i'=1}^{n_m}
    \sum_{j'=1}^{n_\lambda}
    \exp(z_{i',j'})
    },
    \qquad
    i=1,\ldots,n_m,\quad j=1,\ldots,n_\lambda,
    \label{eq:predicted_distribution}
\end{equation}
where $\theta$ denotes the trainable network parameters. 
This architecture allows the model to infer suitable reconstruction parameters from the geometric and oscillatory features of the input profile.
The detailed network architecture is illustrated in Figure~\ref{fig_network_architecture}.

\paragraph{Physics-guided training objective}

Training seeks to minimize a loss function consisting of two components: a reconstruction-based penalty and a probabilistic matching term.  
For the $s$-th sample $\bigl(v^{(s)},u_{\mathrm{ref}}^{(s)}\bigr)$, the first term enforces reconstruction fidelity with
\begin{equation}
    \label{eq:loss_rec}
    \mathcal L_{\mathrm{rec}}^{(s)}(\theta)
    =
    \left\|
        \widetilde v_\theta^{(s)}
        -
        u_{\mathrm{ref}}^{(s)}
    \right\|_2^2,
\end{equation}
with probability-weighted reconstruction defined as
\begin{equation}
    \widetilde v_\theta^{(s)}
    =
    \mathcal R
    \bigl(
        P_\theta(v^{(s)});v^{(s)}
    \bigr)
    =
    \sum_{i=1}^{n_m}
    \sum_{j=1}^{n_\lambda}
    p_\theta(m_i,\lambda_j\mid v^{(s)})
    v_{i,j}^{(s)}.
\end{equation}
The second term matches the predicted distribution to the soft-target distribution through the Kullback-Leibler divergence~\cite{kullback1951information},
\begin{equation}
    \mathcal L_{\mathrm{KL}}^{(s)}(\theta)
    =
    D_{\mathrm{KL}}
    \left(
        P_\gamma^{\star,(s)}
        \,\middle\|\,
        P_\theta(v^{(s)})
    \right)
    \nonumber
    =
    \sum_{i=1}^{n_m}
    \sum_{j=1}^{n_\lambda}
    P_{\gamma,i,j}^{\star,(s)}
    \log
    \frac{
        P_{\gamma,i,j}^{\star,(s)}
    }{
        p_\theta(m_i,\lambda_j\mid v^{(s)})
    }.
\end{equation}
Thus the full loss for the dataset $\mathcal D = \left\{ \bigl(v^{(s)},u_{\mathrm{ref}}^{(s)}\bigr) \right\}_{s=1}^{N}$ is defined as
\begin{equation}
    \mathcal J_{\mathcal D}(\theta)
    =
    \frac{1}{N}
    \sum_{s=1}^{N}
    \left[
        (1-\alpha)
        \mathcal L_{\mathrm{rec}}^{(s)}(\theta)
        +
        \alpha
        \mathcal L_{\mathrm{KL}}^{(s)}(\theta)
    \right],
    \qquad
    0<\alpha<1.
    \label{eq:total_loss}
\end{equation}
The reconstruction term ensures that the probability-weighted output remains close to the reference PDE solution, while the KL term aligns the predicted parameter distribution with the error-induced probability landscape. 
In this sense, the objective is reconstruction-informed and physics-guided: it uses the analytic structure of Gegenbauer reconstruction together with reference solution data to learn physically meaningful parameter distributions.

\subsubsection{Two-stage Pre-training and Fine-tuning Strategy\label{pretraining_and_finetuning}}

To balance reconstruction accuracy and computational cost, we adopt a two-stage training strategy consisting of general pre-training followed by problem-specific fine-tuning, where pre-training is to learn transferable oscillation-aware features and a broad mapping from input profiles to Gegenbauer distributions, while fine-tuning adapts for a specific PDE and data-driven solver. 
The procedure is applied separately to $r\in\{\mathrm{DIOs},\mathrm{SSOs}\}$ since two oscillation classes exhibit different spatial localization and reconstruction-error distributions.

\paragraph{Pre-training stage.}

Let $\mathcal D_{\mathrm{pre}}^{r}$ be the pre-training dataset with corresponding oscillation class $r\in \{\mathrm{DIOs},\mathrm{SSOs}\}$, the pre-trained parameters are then determined by
\begin{equation}
    \theta_{\mathrm{pre}}^{r}
    \in
    \operatorname*{argmin}_{\theta}
    \mathcal J_{\mathcal D_{\mathrm{pre}}^{r}}(\theta),
    \label{eq:pretraining_problem}
\end{equation}
where $\mathcal J_{\mathcal D}$ is the empirical objective in Equation~\eqref{eq:total_loss}. The DIOs and SSOs predictors use the same parametric family and objective, but are optimized over distinct
training distributions.

To distinguish the transferable and problem-dependent components,
we write the trainable parameters during pre-training as
\begin{equation*}
    \theta_{\mathrm{pre}}^{r}
    =
    \bigl(
        \eta_{\mathrm{pre}}^{r},
        \omega_{\mathrm{pre}}^{r}
    \bigr),
\end{equation*}
where $\eta_{\mathrm{pre}}^{r}$ parameterizes the multi-scale feature map and $\omega_{\mathrm{pre}}^{r}$ parameterizes the final map from the extracted features to the logits over $\Gamma$.

\paragraph{Fine-tuning stage.}
Although the pre-training stage equips the model with broad oscillation recognition capabilities, the projected data do not fully characterize the dynamics of a specific target problem. 
Moreover, different PDEs and data-driven solvers may produce distinct oscillation amplitudes, discontinuity geometries, and parameter-error landscapes. 
Fine-tuning allows the model to specialize to these problem-dependent features while preserving the robust feature extraction learned during pre-training~\cite{wang2026pretraining}.

Let $\mathcal D_{\mathrm{tar}}^{r}$ be the problem-specific dataset for a target PDE and data-driven solver.
In the lightweight strategy adopted here, the feature parameters are fixed at $\eta_{\mathrm{pre}}^{r}$ and only the final parameter map is optimized:
\begin{equation}
    \omega_{\mathrm{ft}}^{r}
    \in
    \operatorname*{argmin}_{\omega}
    \mathcal J_{\mathcal D_{\mathrm{tar}}^{r}}
    \bigl(
        \eta_{\mathrm{pre}}^{r},\omega
    \bigr),
    \label{eq:finetuning_problem}
\end{equation}
and the resulting fine-tuned parameters are
\begin{equation*}
    \theta_{\mathrm{ft}}^{r}
    =
    \bigl(
        \eta_{\mathrm{pre}}^{r},
        \omega_{\mathrm{ft}}^{r}
    \bigr).
\end{equation*}
Thus, fine-tuning preserves the generic multi-scale representation learned during pre-training and recalibrates only the map from the extracted features to the probability distribution over Gegenbauer parameters.

\vspace{1.0em}

Overall, the two-stage strategy provides a practical compromise between full re-training and direct use of the pre-trained model. 
It retains the general feature representation learned from large-scale pre-training data, while using a small amount of problem-specific data to improve reconstruction accuracy and reduce computational cost.

    \section{Numerical studies \label{numerical_results}}

This section evaluates the proposed probabilistic Gegenbauer reconstruction framework on one- and two-dimensional transport-dominated problems. 
Subsection~\ref{General_settigs} details the pre-training procedure. 
Subsection~\ref{test_problmes} assesses the reconstruction performance on several representative test problems, including both linear and nonlinear cases. 
Subsection~\ref{Accuracy_and_efficiency_tradeoff} further investigates the cost-accuracy trade-off, fine-tuning strategies, and robustness of the proposed two-stage training framework.

To quantitatively assess reconstruction quality, the relative $\ell^2$ and maximum errors are computed on a masked index set $\mathcal I$, where five grid points on each side of each detected discontinuity are excluded to reduce the influence of small localization mismatches, 
\begin{equation*}
    e_{\mathrm{rel}} = \frac{\left(\sum_{i\in\mathcal I} |u_i^{\mathrm{recon}}-u_i^{\mathrm{true}}|^2\right)^{1/2}}{\left(\sum_{i\in\mathcal I}|u_i^{\mathrm{true}}|^2\right)^{1/2}},
    \qquad
    e_{\max} = \max_{i\in\mathcal I}|u_i^{\mathrm{recon}} - u_i^{\mathrm{true}}|,
\end{equation*}
where $u^{\text{recon}}_{i}$ and $u^{\text{true}}_{i}$ denote the $i$-th component of the reconstructed and reference solutions, respectively, with the same spatial mask applied.

To further evaluate the generalization performance over multiple test samples, we report statistical summaries of the following metrics:
\begin{eqnarray*}
    && \text{Mean} = \frac{1}{N}\sum_{j=1}^{N} e_j, \\[4pt]
    && \text{Std} = \left( \frac{1}{N}\sum_{j=1}^{N} (e_j - \text{Mean})^2 \right)^{1/2}, \\[4pt]
    && \text{5\% quantile} = \inf\left\{\epsilon : \frac{ \#\{j : e_j \le \epsilon\}}{N} \ge 0.05 \right\}, \\[4pt]
    && \text{95\% quantile} = \inf\left\{\epsilon : \frac{\#\{j : e_j \le \epsilon\}}{N} \ge 0.95 \right\},
\end{eqnarray*}
where $\{e_j\}_{j=1}^{N}$ is the given set of error values computed over $N$ samples.

Experiments are conducted on a Linux server running Ubuntu~20.04 (kernel 5.15.0-67-generic) with Python~3.12.7 (Anaconda). The computational environment includes PyTorch~2.6.0 with CUDA~12.4, NumPy~1.26.4, and SciPy~1.13.1. The server is equipped with an NVIDIA RTX~A6000 GPU and 512~GB of RAM.

\subsection{General Settings and Pre-training Performance\label{General_settigs}}

This subsection details the common experimental settings and reports the performance of the pre-trained parameter predictors. 

All training datasets are generated on a uniform grid with $n_{\mathrm{space}}=256$ points over $x\in[-1,1)$. 
The candidate Gegenbauer parameter set is spanned with $m=0,\ldots,10$ and $\lambda=1,\ldots,20$.
To separate discontinuity-induced and spurious oscillatory samples, we use the oscillation-discrimination criterion described in Section~\ref{oscillation_discrimination}, with $S_{\mathrm{thresh}}=3.0$ and $\mathrm{edge}_{\mathrm{frac}}=0.25$. 
The detailed construction of the discontinuity-induced and spurious pre-training datasets, including the sampling procedures and parameter choices, is provided in Appendix~\ref{appendix_dataset_pretraining}.
The total dataset-generation time is $370.42\,\mathrm{s}$ for the discontinuity-induced dataset and $374.57\,\mathrm{s}$ for the spurious dataset during pre-training.

Two pre-trained predictors are constructed, corresponding to discontinuity-induced oscillations and spurious numerical oscillations. Both predictors use the same multi-scale CNN architecture, \texttt{CNN\_ReconParamNet}, introduced in Section~\ref{prob_model_and_loss} and detailed in Appendix~\ref{pretrain_network}. 
All pre-training, fine-tuning, and re-training experiments use the same network architecture.

Training is performed with a batch size of $32$ using the AdamW optimizer and a cosine-annealing learning-rate schedule, with the learning rate decreasing from $10^{-3}$ to $10^{-5}$ over $600$ epochs. 
The total training time is approximately $308.37\,\mathrm{s}$ for the discontinuity-induced dataset and $307.68\,\mathrm{s}$ for the spurious dataset.
For both oscillation types, the dataset is split to 80\% for training and 20\% for validation, and the checkpoint with the lowest validation loss is selected as the pre-trained model.

To assess the generalization performance of the pre-trained models, $100$ additional test samples are generated for each oscillation type. 
The corresponding error statistics are summarized in Table~\ref{general_model_error_statistics}.

\begin{table}[!ht]
    \centering
    \small
    \caption{\label{general_model_error_statistics}Error statistics of the pre-trained models evaluated on $100$ newly generated test samples for each oscillation type. 
    }
    \begin{tabular}{c|l|c|c|c|c}
        \hline
        & & Mean & Std & 5\% quantile & 95\% quantile\\\hline
        \multirow{2}{*}{DIOs} 
        & Relative $\ell^2$& 1.594850e-02&1.425302e-02&2.669373e-03&4.686053e-02\\
        & Maximum &2.542741e-01&2.983944e-01&4.198718e-02&7.650095e-01\\\hline
        \multirow{2}{*}{SOSs} 
        & Relative $\ell^2$ & 9.199257e-03&8.445588e-03&1.372313e-03&2.582346e-02 \\
        & Maximum &5.208429e-03&3.335678e-03&8.909389e-04&1.072704e-02\\\hline
    \end{tabular}   
\end{table}

\subsection{Test Problems and Qualitative Results\label{test_problmes}}

We evaluate the proposed physics-informed post-processing framework on representative transport-dominated benchmarks, including linear and nonlinear problems. 
The purpose of these experiments is not to develop new data-driven solvers, but to assess whether the proposed probabilistic Gegenbauer reconstruction can effectively improve oscillatory data-driven solutions across different model classes and problem settings.
For all test problems, we report results for G-ROM solutions, which serve as representative reduced-order approximations. 
To further demonstrate applicability beyond reduced-order modeling, we also apply the proposed correction to DeepONet predictions for the one-dimensional Burgers' equation. 
For all G-ROM examples, snapshots are generated from exact solutions when available; otherwise, the numerical solutions are obtained by the fifth-order WENO finite difference scheme in space and third-order Runge-Kutta method in time.
We remark that if we use a numerical method with damping to obtain the solution snapshots, the resulting basis functions are qualitatively similar to those obtained by exact solutions with the pre-processing filter.

\subsubsection{\texorpdfstring{ Example 1: linear equation}{Example 1: linear equation}}\label{example1}

We first consider a parametric linear transport equation,
\begin{equation}
    \label{eg1_hyperbolic_equation}
    \frac{\partial u}{\partial t} = -2 \frac{\partial u}{\partial x}, \quad x\in[-1,1],
\end{equation}
with initial data $u(x,0)$ given by random polynomials of degree $\le 10$.
The exact solution of Equation~\eqref{eg1_hyperbolic_equation} can be expressed as the shift profile $ u(x,t) = u(x - 2 t, 0) $.
This example provides a simple benchmark for assessing the proposed post-processing framework on G-ROM solutions with transport-induced oscillations.

A POD-Galerkin reduced-order model is used as the oscillatory baseline. 
The reduced space is constructed from filtered solution snapshots, and the G-ROM solution at target time $T=0.2$ serves as the input to the proposed post-processing framework. 
The detailed G-ROM construction is provided in Appendix~\ref{appendix_eg1_grom} and the problem-specific fine-tuning and re-training adaptation setups are shown in Appendix~\ref{appendix_eg1_adaptation_setup}.

\begin{figure}[!ht]
    \centering
    \subfloat[Empirical Gegenbauer sweep]{\includegraphics[width=0.31\textwidth]{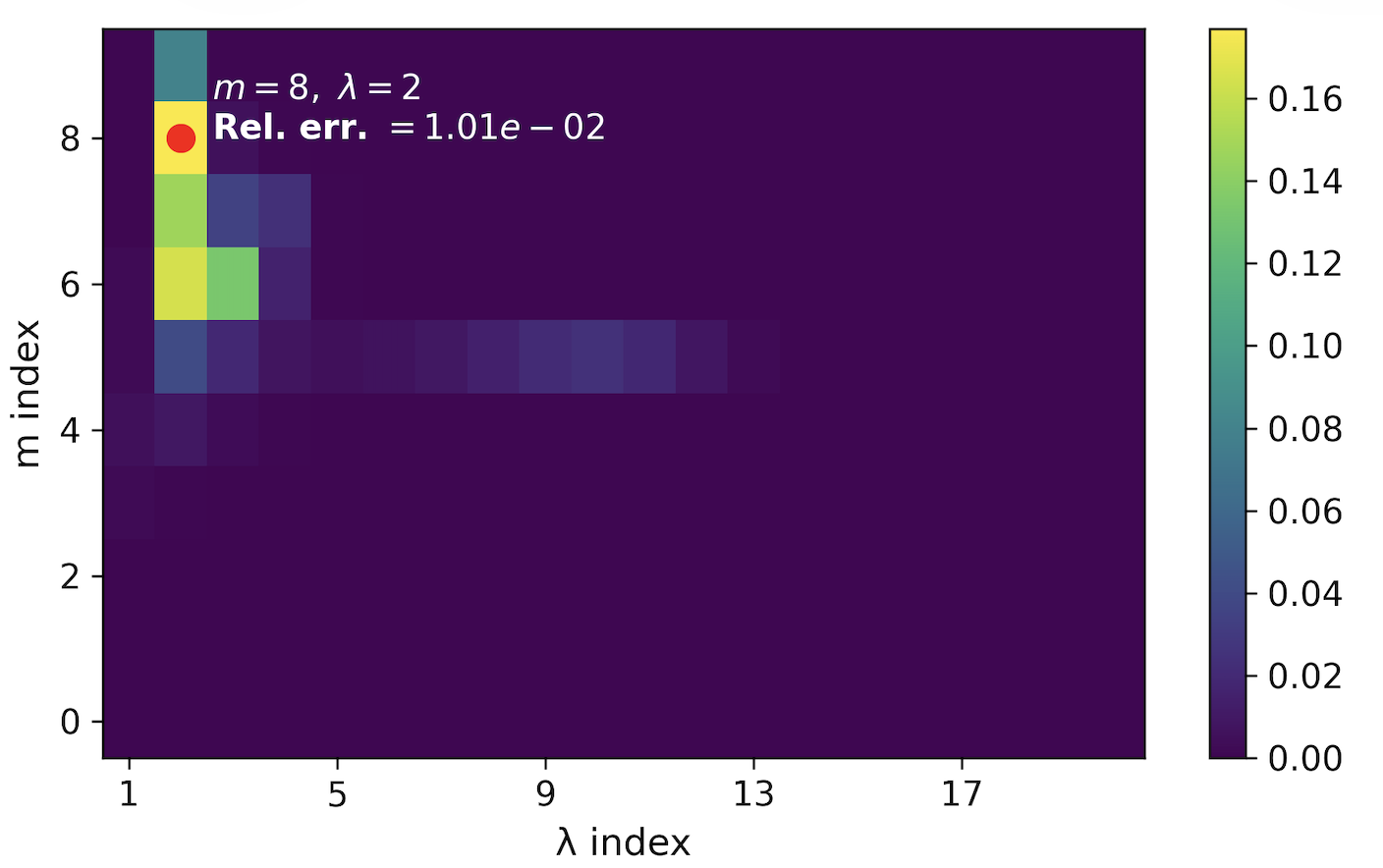}}
    \hspace{0.02\textwidth}
    \subfloat[Empirical best reconstruction]{\includegraphics[width=0.32\textwidth]{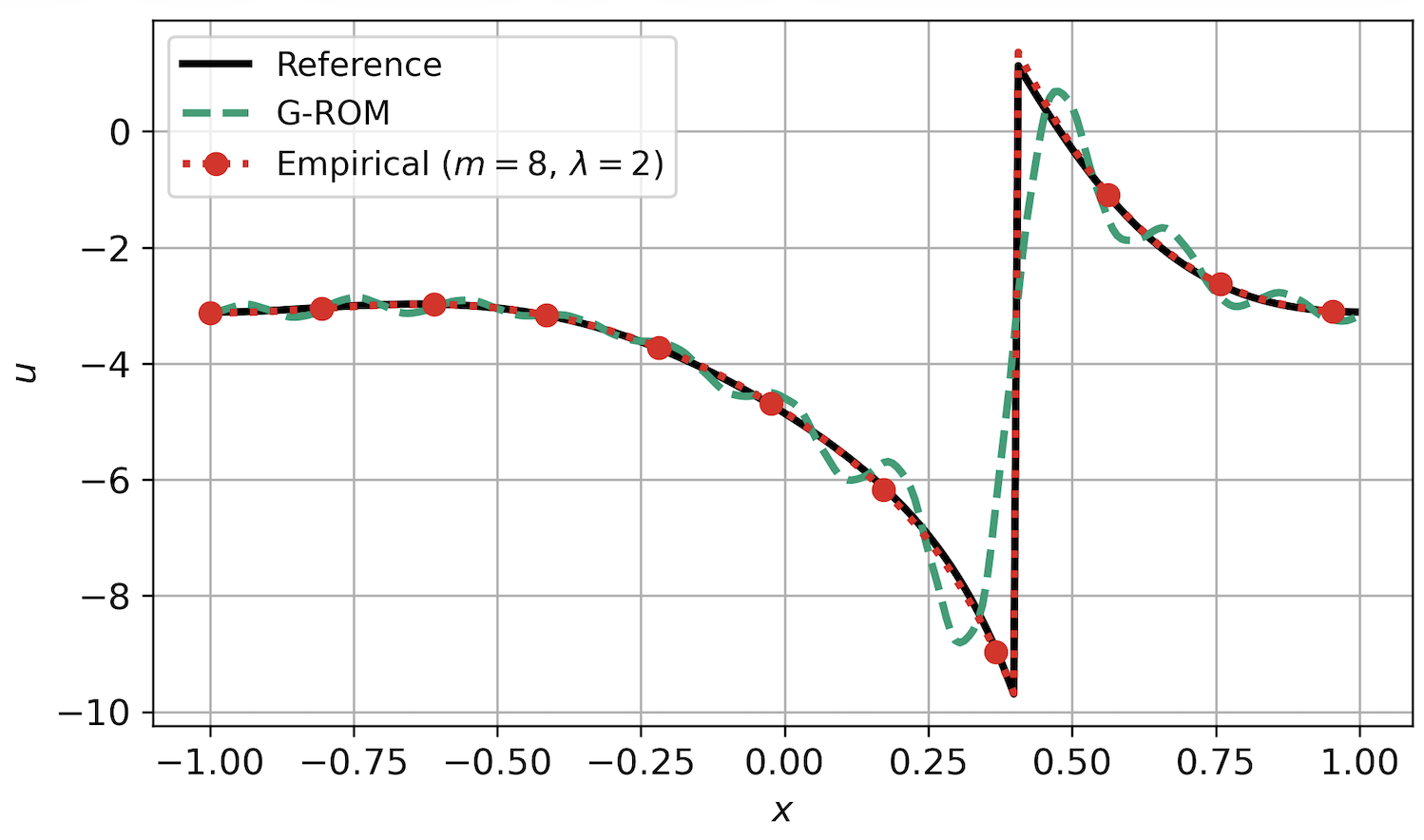}}
    \\[2mm]
    \subfloat[Pre-trained model]{
        \includegraphics[width=0.3\textwidth]{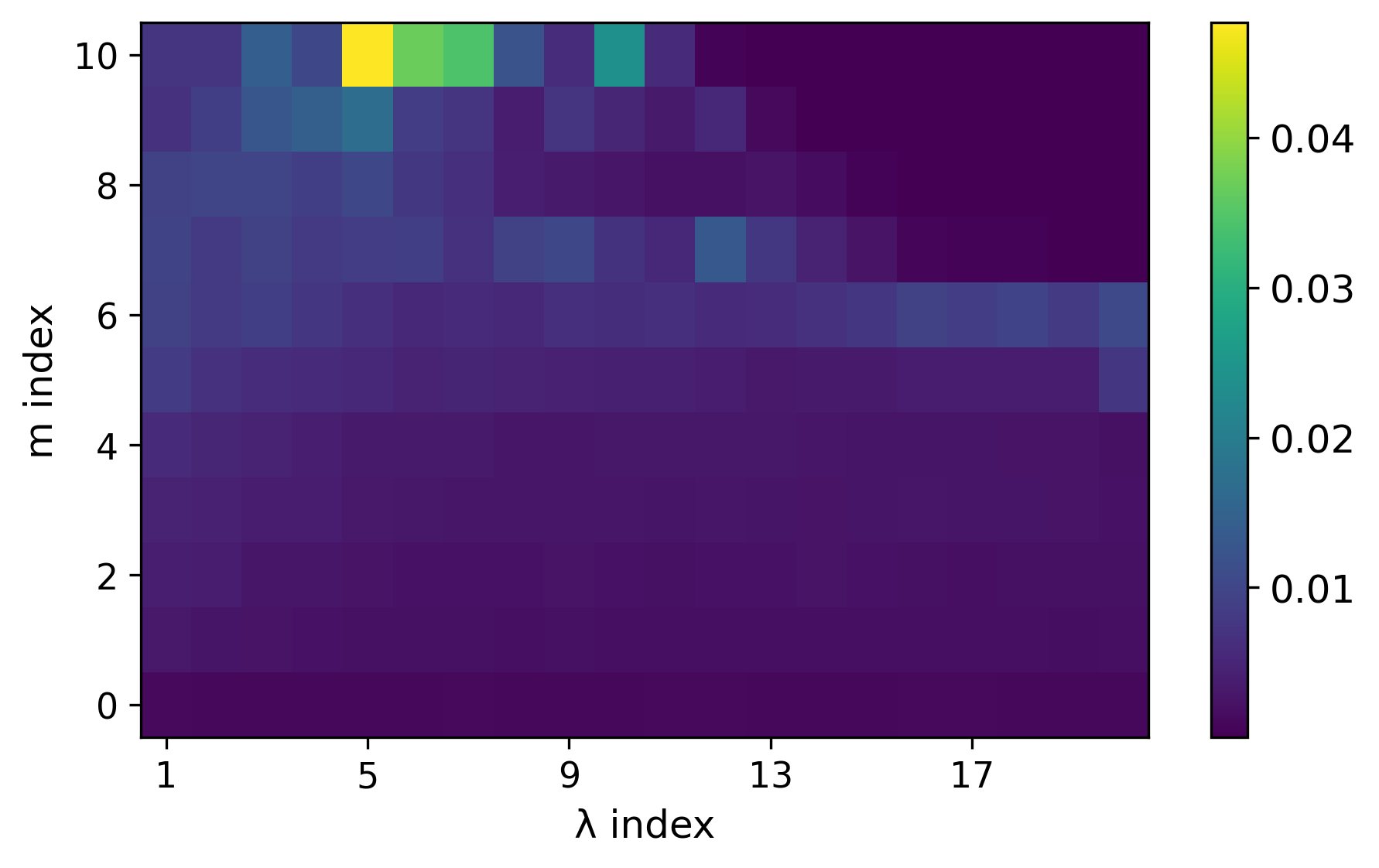}
    }
    \hspace{0.01\textwidth}
    \subfloat[Fine-tuned model]{
        \includegraphics[width=0.3\textwidth]{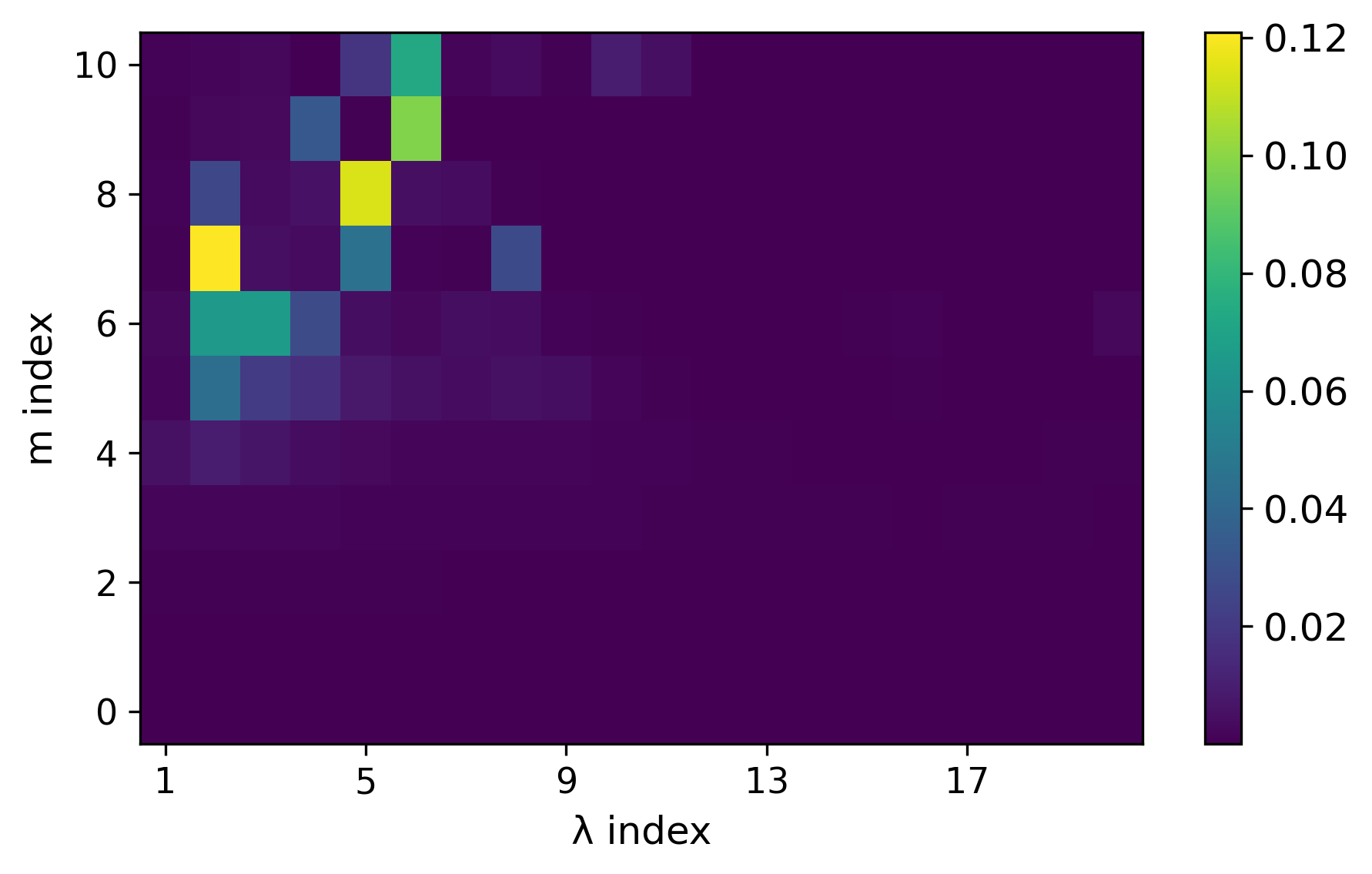}
    }
    \hspace{0.01\textwidth}
    \subfloat[Re-trained model]{
        \includegraphics[width=0.3\textwidth]{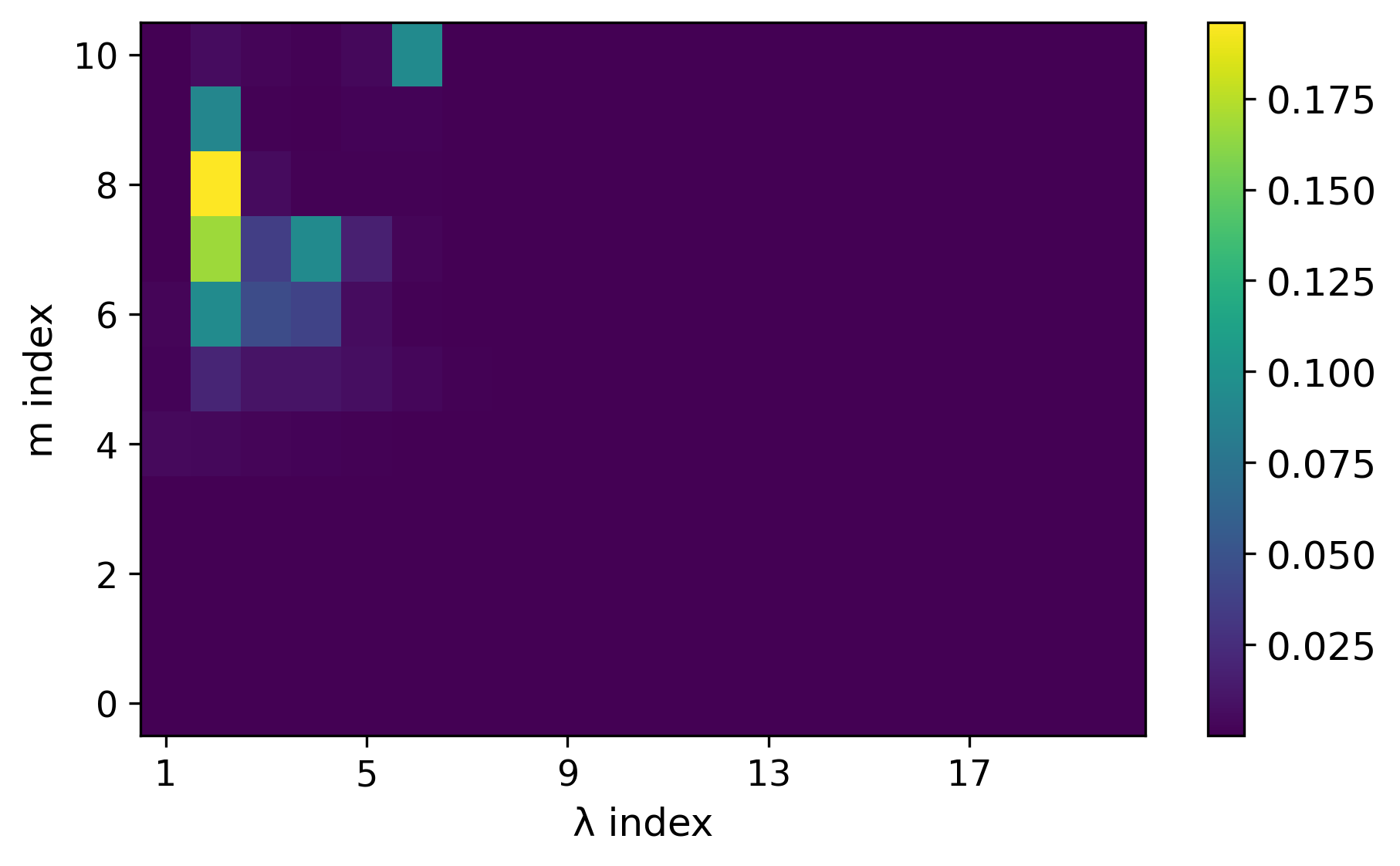}
    }
    \caption{Example 1. Parameter-distribution diagnostics for the representative test sample: empirical Gegenbauer sweep with best reconstruction, and predicted distributions from the pre-trained, fine-tuned, and re-trained models. }
    \label{eg1_hyperbolic_parameter_diagnostics}
\end{figure}

\begin{figure}[!ht]
    \centering
    \subfloat[Pre-trained model]{\includegraphics[width=0.3\textwidth]{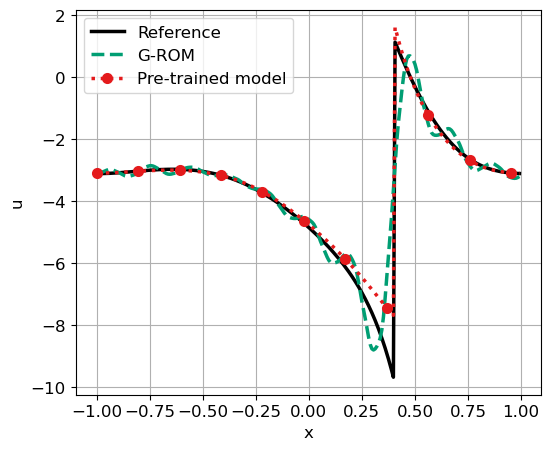}}
    \hspace{0.03\textwidth}
    \subfloat[Fine-tuned model]{\includegraphics[width=0.3\textwidth]{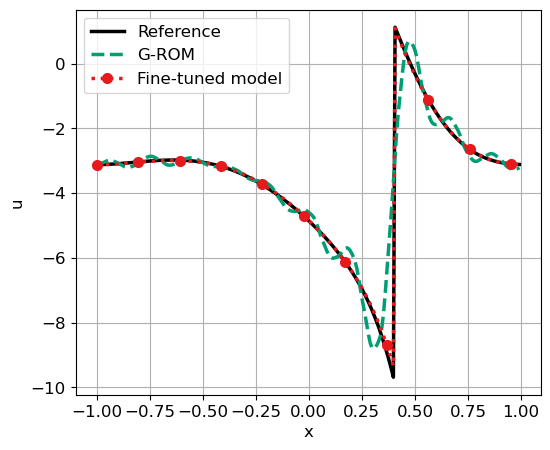}}
    \hspace{0.03\textwidth}
    \subfloat[Re-trained model]{\includegraphics[width=0.3\textwidth]{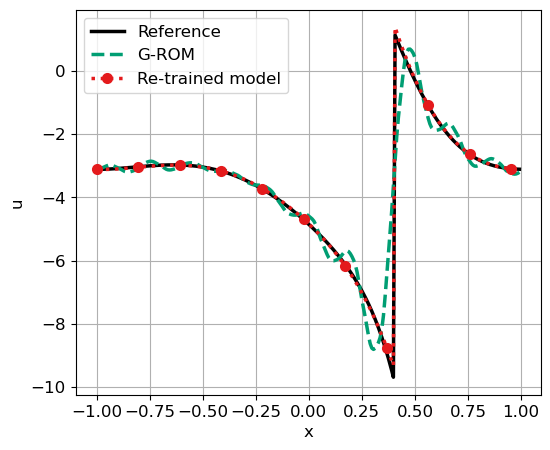}}
    \caption{Example 1. Reconstructions produced by the pre-trained, fine-tuned, and re-trained models for a randomly selected test sample.}
    \label{eg1_hyperbolic_specific_sample}
\end{figure}

To illustrate the reconstruction behavior, we first consider a representative G-ROM solution at $T=0.2$ with polynomial initial condition of $\mathrm{degree}=8$. 
Figure~\ref{eg1_hyperbolic_parameter_diagnostics} provides the parameter-distribution diagnostics for this representative sample. 
The empirical Gegenbauer sweep indicates that several parameter pairs near the optimum yield comparable reconstruction quality, motivating a weighted reconstruction rather than hard best-pair selection. 
Compared with the pre-trained model, the fine-tuned and re-trained models produce distributions that better match the empirical landscape, showing that problem-specific adaptation effectively captures the target reconstruction behavior and Figure~\ref{eg1_hyperbolic_specific_sample} shows the corresponding probability-weighted reconstructions, with quantitative error metrics reported in Table~\ref{eg1_hyperbolic_specific_sample_error}. 
All three models reduce the oscillations in the raw G-ROM solution, while the fine-tuned and re-trained models provide sharper reconstructions and closer agreement with the reference solution.

\begin{table}[!ht]
    \centering
    \caption{\label{eg1_hyperbolic_specific_sample_error} Example 1. Relative $\ell^2$ and maximum errors of the G-ROM solution and the reconstructions obtained from the pre-trained, fine-tuned, re-trained models and best empirical reconstruction for the randomly selected test sample.}
    \begin{tabular}{c|c|c|c|c|c}\hline
        & G-ROM & Pre-trained  & Fine-tuned & Re-trained & Empirical \\
        \hline
        Relative $\ell^2$& 1.2267e-01  &7.3272e-02  &1.0439e-02 & 9.5672e-03 & 8.7081e-03 \\
        \hline
        Maximum&4.1858e+00 &1.7149e+00 & 2.8026e-01 & 2.5848e-01 & 1.5940e-01 \\
        \hline
    \end{tabular}
\end{table}

\begin{figure}[!ht]
    \centering
    \includegraphics[width=0.9\linewidth]{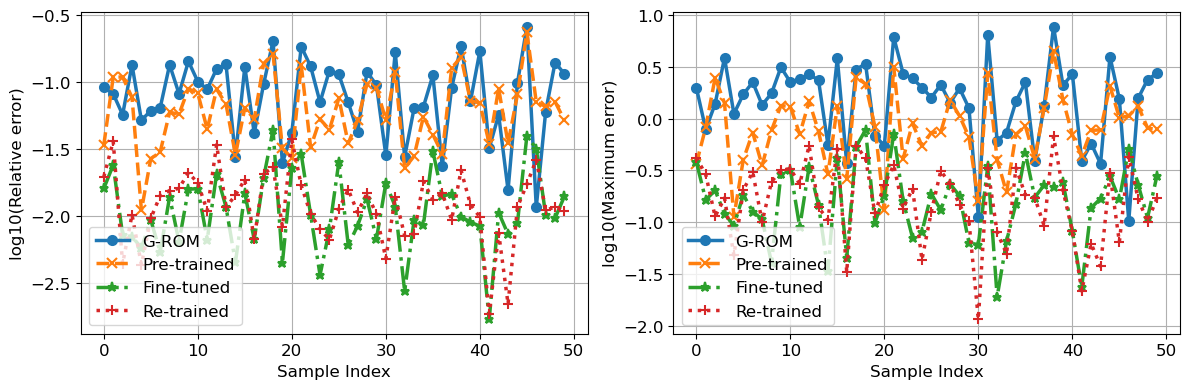}
    \caption{Example 1. 
    Per-sample errors over $50$ test samples for the pre-trained, fine-tuned, and re-trained models. Left: relative $\ell^2$ errors; right: maximum errors}
    \label{eg1_hyperbolic_model_comparision}
\end{figure}

To assess statistical performance, we evaluate the models on $50$ randomly generated test samples at $T=0.2$. Table~\ref{eg1_hyperbolic_model_comparision_errors} reports the mean, standard deviation, and empirical $5\%$ and $95\%$ quantiles of the relative $\ell^2$ and maximum errors. 
The corresponding per-sample error curves are provided in Figure~\ref{eg1_hyperbolic_model_comparision}. 
The fine-tuned and re-trained models substantially reduce both error metrics compared with the raw G-ROM and pre-trained reconstructions. Moreover, their performance is nearly indistinguishable, indicating that fine-tuning is sufficient to adapt the pre-trained predictor to this linear transport regime.

\begin{table}[!ht]
    \centering
    \caption{\label{eg1_hyperbolic_model_comparision_errors} Example 1. 
    Statistical summary of the relative $\ell^2$ and maximum errors for the G-ROM solutions and the reconstructed solutions over $50$ test samples.
    }
    \begin{tabular}{c|l|c|c|c|c}
        \hline
        & & Mean & Std & 5\% quantile & 95\% quantile\\\hline
        \multirow{2}{*}{G-ROM} 
        & Relative $\ell^2$& 9.366e-02& 5.356e-02& 2.427e-02&1.860e-01\\
        & Maximum &2.054e+00&1.564e+00&3.682e-01& 5.161e+00\\\hline
        \multirow{2}{*}{Pre-trained model} 
        & Relative $\ell^2$ &  7.044e-02&4.293e-02& 2.698e-02&1.470e-01 \\
        & Maximum &1.056e+00&8.474e-01&1.751e-01&2.640e+00\\\hline
        \multirow{2}{*}{Fine-tuned model} 
        & Relative $\ell^2$ &  1.327e-02& 9.168e-03 & 4.003e-03&3.132e-02 \\
        & Maximum &2.099e-01&1.652e-01&3.701e-02& 5.258e-01\\\hline
        \multirow{2}{*}{Re-trained model} 
        & Relative $\ell^2$ & 1.385e-02&7.735e-03 & 4.343e-03&3.065e-02 \\
        & Maximum &2.136e-01&1.540e-01&3.503e-02&5.338e-01\\\hline
    \end{tabular}   
\end{table}

Table~\ref{eg1_hyperbolic_time_tracing} compares the computational costs of fine-tuning and full re-training. 
Although both approaches achieve comparable reconstruction accuracy, fine-tuning uses only $100$ samples and requires substantially less training time than re-training on $2000$ samples. This confirms that the proposed two-stage strategy provides a favorable cost-accuracy trade-off for this benchmark.

\begin{table}[!ht]
    \centering
    \caption{\label{eg1_hyperbolic_time_tracing} Example 1. 
    Comparison of dataset generation and training costs for fine-tuned and re-trained models in the G-ROM setting for 1D hyperbolic equation. }
    \begin{tabular}{c|c|c}\hline
         & Dataset generation & Model training  \\ \hline
         Fine-tuned model & 22.80 s / 100 samples & 7.71 s / 200 epochs\\ \hline
         Re-trained model & 406.90 s / 2000 samples & 296.63 s / 600 epochs\\\hline
    \end{tabular}
\end{table}

\subsubsection{Example 2: inviscid Burgers' equation\label{eg2_1d_burgers}}

We next consider the inviscid Burgers' equation with moving discontinuity location:
\begin{equation}
    \begin{cases}
    u_t + uu_x = 0, \qquad x \in [0,1],\\[2mm]
    u(x,0) = a + b\sin(\omega x + \varphi), 
    \end{cases}
\end{equation}
subjected to periodic boundary conditions, with the parameters chosen as $a \in (-2,2),\; b \in (0.1,2),\; \varphi \in (-\pi,\pi),\; \omega = \pi$.
We consider two representative baselines: a POD-Galerkin ROM and a DeepONet surrogate. 
Implementation details for both baselines and their problem-specific adaptation datasets are provided in Appendix~\ref{appendix_eg2}.
The G-ROM solution at $T=0.6$, for which the selected samples have already formed shocks, is used as the oscillatory input to the proposed post-processing framework. 
Details of the G-ROM construction, fine-tuning and re-training setups
are provided in Appendix~\ref{appendix_eg2_grom}.

Figure~\ref{eg2_burgers_specific_sample} and Table~\ref{eg2_burgers_specific_sample_error} show the reconstruction performance for a representative test sample. 
The pre-trained model exhibits noticeable reconstruction errors, whereas both the fine-tuned and re-trained models provide highly accurate reconstructions, with the re-trained model achieving the smallest relative $\ell^2$ and maximum errors.
This indicates that problem-specific adaptation is necessary for the Burgers setting.

\begin{figure}[!htbp]
    \centering
    \subfloat[Pre-trained model]{\includegraphics[width=0.3\textwidth]{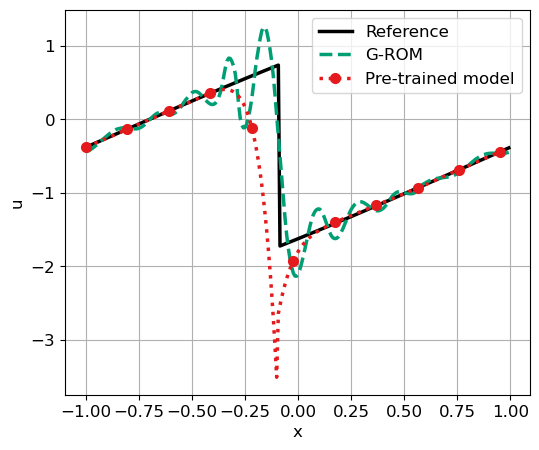}}
    \hspace{0.03\textwidth}
    \subfloat[Fine-tuned model]{\includegraphics[width=0.3\textwidth]{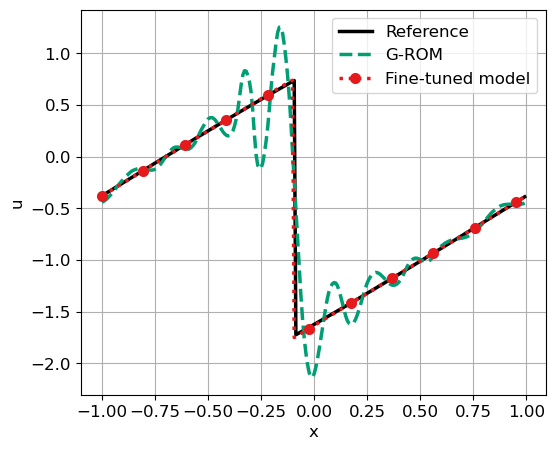}}
    \hspace{0.03\textwidth}
    \subfloat[Re-trained model]{\includegraphics[width=0.3\textwidth]{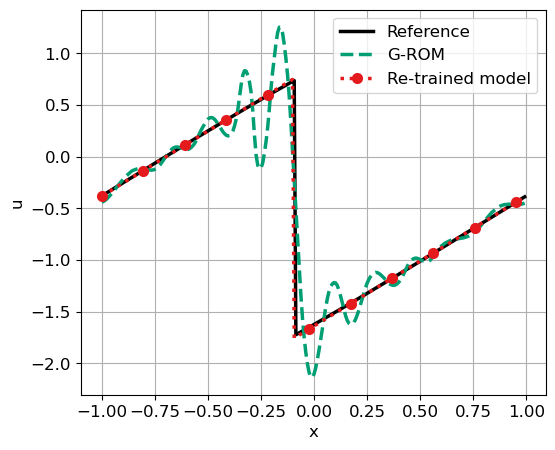}}
    \caption{Example 2. Reconstructions of G-ROM solutions produced by the pre-trained, fine-tuned, and re-trained models for a selected test sample.}
    \label{eg2_burgers_specific_sample}
\end{figure}

\begin{table}[!htbp]
    \centering
    \caption{\label{eg2_burgers_specific_sample_error} Example 2. 
    Relative $\ell^2$ and maximum errors of the G-ROM solutions and the reconstructions for the selected test sample.}
    \begin{tabular}{c|c|c|c|c}\hline
        & G-ROM & Pre-trained model & Fine-tuned model & Re-trained model \\
        \hline
        Relative $\ell^2$& 2.1469e-01  & 5.5649e-01  &5.2751e-03 &  6.7094e-03  \\
        \hline
        Maximum&6.6364e-01 & 3.4263e+00 &1.5862e-02 &2.0901e-02\\
        \hline
    \end{tabular}
\end{table}

\begin{figure}[!ht]
    \centering
    \includegraphics[width=0.9\linewidth]{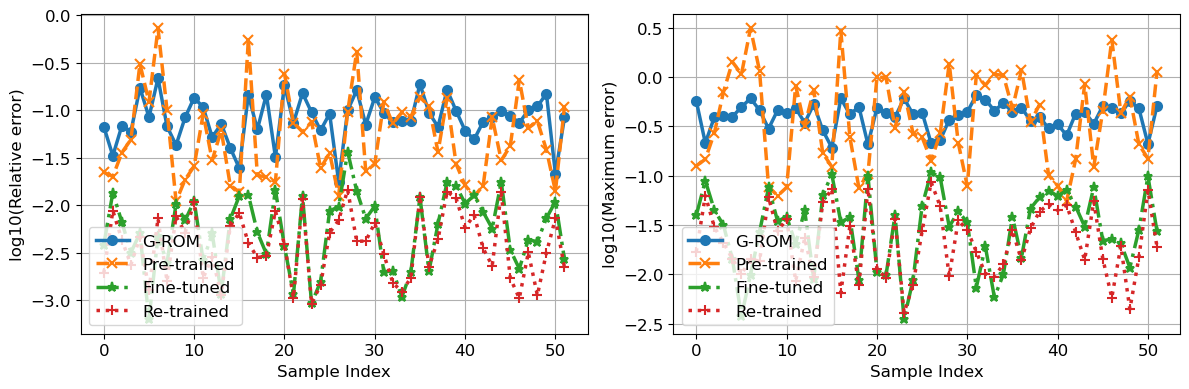}
    \caption{Example 2. 
    Per-sample errors over $52$ G-ROM test samples for the pre-trained, fine-tuned, and re-trained models.  Left: relative $\ell^2$ errors; right: maximum errors.
    }
    \label{eg2_burgers_model_comparision}
\end{figure}

To further evaluate generalization performance, $52$ additional discontinuous test samples are generated and used to assess all three models. 
As summarized in Table~\ref{eg2_burgers_model_comparision_errors}, both the fine-tuned and re-trained models consistently outperform the pre-trained model across all error metrics. 
Specifically, the fine-tuned model reduces the mean relative $\ell^2$ error from $9.197\times10^{-2}$ for the raw G-ROM solution to $7.508\times10^{-3}$, and achieves accuracy close to the fully re-trained model. 
The corresponding per-sample error curves are provided in Figure~\ref{eg2_burgers_model_comparision}.

\begin{table}[!ht]
    \centering
    \caption{\label{eg2_burgers_model_comparision_errors} Example 2. 
    Statistical summary of the relative $\ell^2$ and maximum errors for the G-ROM solutions and the reconstructed solutions over 52 test samples.}
    \begin{tabular}{c|l|c|c|c|c}
        \hline
        & & Mean & Std & 5\% quantile & 95\% quantile\\\hline
        \multirow{2}{*}{G-ROM} 
        & Relative $\ell^2$& 9.197e-02& 4.523e-02& 2.871e-02& 1.769e-01\\
        & Maximum &4.230e-01&1.149e-01&2.116e-01&6.079e-01\\\hline
        \multirow{2}{*}{Pre-trained model} 
        & Relative $\ell^2$ & 9.367e-02&1.339e-01&1.326e-02&3.548e-01 \\
        & Maximum &6.286e-01&6.795e-01&7.577e-02&1.863e+00\\\hline
        \multirow{2}{*}{Fine-tuned model} 
        & Relative $\ell^2$ &7.508e-03&6.247e-03&1.131e-03&1.668e-02 \\
        & Maximum &3.955e-02&2.898e-02&6.655e-03&9.943e-02\\\hline
        \multirow{2}{*}{Re-trained model} 
        & Relative $\ell^2$ & 5.012e-03&3.783e-03& 1.106e-03&1.293e-02\\
        & Maximum &2.793e-02&2.096e-02&6.170e-03&7.210e-02\\\hline
    \end{tabular}   
\end{table}

Beyond reduced-order modeling, we further apply the proposed probabilistic Gegenbauer reconstruction to solutions generated by a DeepONet surrogate. The DeepONet approximates the solution operator mapping the initial condition to the solution at the target time $T=0.6$, and the detailed architecture, training procedure, and problem-specific adaptation setup are provided in Appendix~\ref{appendix_eg2_deeponet}.

\begin{figure}[!ht]
    \centering
    \subfloat[Pre-trained model]{\includegraphics[width=0.3\textwidth]{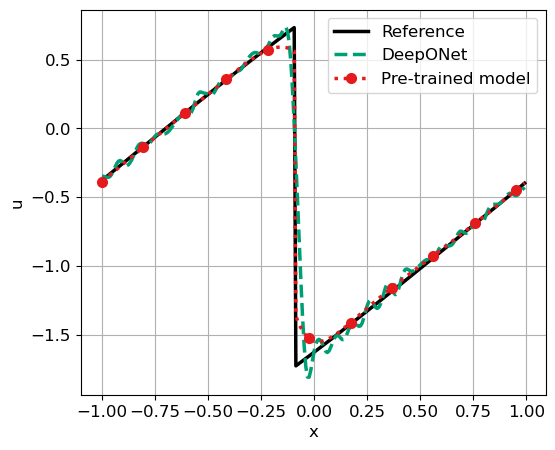}}
    \hspace{0.03\textwidth}
    \subfloat[Fine-tuned model]{\includegraphics[width=0.3\textwidth]{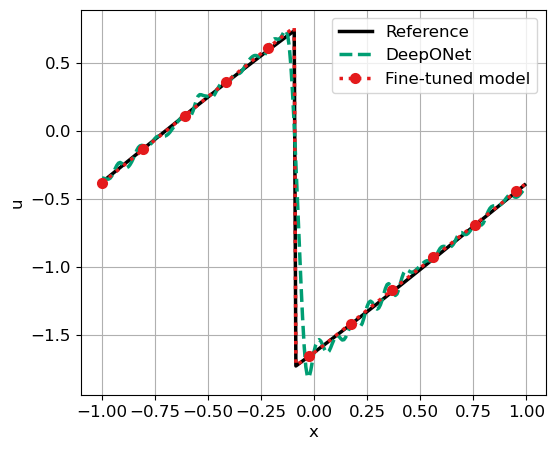}}
    \hspace{0.03\textwidth}
    \subfloat[Re-trained model]{\includegraphics[width=0.3\textwidth]{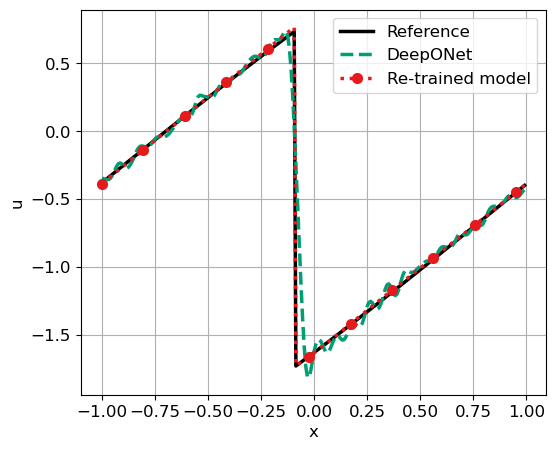}}
    \caption{Example 2. Reconstructions of DeepONet solutions produced by the pre-trained, fine-tuned, and re-trained models for a randomly selected test sample.}
    \label{eg2_burgers_specific_sample_deeponet}
\end{figure}

\begin{figure}[!ht]
    \centering
    \includegraphics[width=0.9\linewidth]{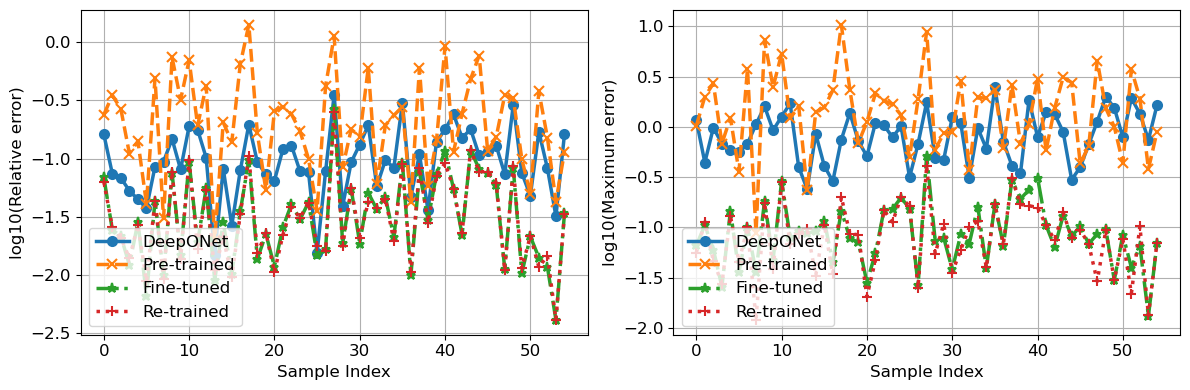}
    \caption{Example 2. 
    Per-sample errors over $55$ DeepONet test samples for the pre-trained, fine-tuned, and re-trained models. Left: relative $\ell^2$ errors; right: maximum errors.
    }
    \label{eg2_burgers_model_comparision_deeponet}
\end{figure}

Representative reconstructions for a selected test case are shown in Figure~\ref{eg2_burgers_specific_sample_deeponet}, with pointwise errors reported in Table~\ref{eg2_burgers_specific_sample_error_deeponet}. 
To evaluate whether the improvement is consistent, we further test the models on $55$ additional DeepONet samples with shock formation, with the statistical results summarized in Table~\ref{eg2_burgers_model_comparision_errors_deeponet}. 
Both the fine-tuned and re-trained models substantially reduce the relative $\ell^2$ and maximum errors compared with the raw DeepONet predictions and the directly pre-trained correction. 
Their mean relative errors are nearly identical, indicating that problem-specific fine-tuning achieves re-training-level accuracy while requiring significantly fewer samples. The corresponding per-sample error curves are provided in Figure~\ref{eg2_burgers_model_comparision_deeponet}.

\begin{table}[!ht]
    \centering
    \caption{\label{eg2_burgers_specific_sample_error_deeponet} Example 2. 
    Relative $\ell^2$ and maximum errors of the DeepONet solutions and the reconstructions for the selected test sample.}
    \begin{tabular}{c|c|c|c|c}\hline
        & DeepONet & Pre-trained model & Fine-tuned model & Re-trained model \\
        \hline
        Relative $\ell^2$& 9.1899e-02 &4.8119e-02  &1.0702e-02 & 1.1305e-02  \\
        \hline
        Maximum&8.9461e-01 & 2.8414e-01 &2.4568e-02 &2.9559e-02\\
        \hline
    \end{tabular}
\end{table}

\begin{table}[!ht]
    \centering
    \caption{\label{eg2_burgers_model_comparision_errors_deeponet} Example 2. 
    Statistical summary of the relative $\ell^2$ and maximum errors for the DeepONet solutions and the reconstructed solutions over $55$ test samples.}
    \begin{tabular}{c|l|c|c|c|c}
        \hline
        & & Mean & Std & 5\% quantile & 95\% quantile\\\hline
        \multirow{2}{*}{G-ROM} 
        & Relative $\ell^2$& 1.121e-01&7.092e-02&3.044e-02&2.548e-01\\
        & Maximum &9.432e-01&5.111e-01&3.094e-01&1.878e+00\\\hline
        \multirow{2}{*}{Pre-trained model} 
        & Relative $\ell^2$ & 2.957e-01&2.843e-01&3.954e-02&8.125e-01 \\
        & Maximum &2.029e+00&2.000e+00&3.690e-01&5.945e+00\\\hline
        \multirow{2}{*}{Fine-tuned model} 
        & Relative $\ell^2$ &4.353e-02&4.264e-02&9.704e-03&1.036e-01 \\
        & Maximum &1.100e-01&8.810e-02&2.776e-02&2.928e-01\\\hline
        \multirow{2}{*}{Re-trained model} 
        & Relative $\ell^2$ &4.260e-02&4.066e-02& 9.516e-03&9.925e-02\\
        & Maximum &1.001e-01&7.586e-02& 2.157e-02&2.261e-01\\\hline
    \end{tabular}   
\end{table}

Table~\ref{eg2_burgers_cost_grom_deeponet} reports the dataset-generation and training costs for fine-tuning and full re-training in both G-ROM and DeepONet settings. 
For both baselines, fine-tuning requires substantially fewer samples and much shorter training time than re-training. 
Together with the error statistics reported above, these results demonstrate that the proposed two-stage strategy achieves a favorable cost--accuracy trade-off for nonlinear Burgers dynamics.

\begin{table}[!ht]
    \centering
    \caption{Example 2. Dataset generation and training costs for fine-tuning and re-training in the G-ROM and DeepONet settings for 1D Burgers' equation.}
    \label{eg2_burgers_cost_grom_deeponet}
    \begin{tabular}{c|c|c|c}
    \hline
    Baseline & Model & Dataset generation & Model training \\
    \hline
    \multirow{2}{*}{G-ROM}
    & Fine-tuned & 57.46 s / 158 samples & 9.81 s / 200 epochs \\
    & Re-trained & 664.15 s / 1926 samples & 286.66 s / 600 epochs \\
    \hline
    \multirow{2}{*}{DeepONet}
    & Fine-tuned & 42.25 s / 158 samples & 8.09 s / 200 epochs \\
    & Re-trained & 384.81 s / 1929 samples & 281.73 s / 600 epochs \\
    \hline
    \end{tabular}
\end{table}

\subsubsection{Example 3: 2D linear equation\label{eg3_hyperbolic_2d}}
We consider the two-dimensional linear hyperbolic equation
\begin{equation}
    \label{eq_hyper_2d}
    u_t + uu_x + uu_y = 0,\qquad (x,y)\in [-1,1]\times [-1,1],
\end{equation}
with the random initial condition given by
\begin{equation}    
    \label{eq_hyper_2d_ic}
    \begin{gathered}
        u(x,0) = \begin{cases}
            a + b \sin(\omega_1 x)\cos(\omega_2 y),\qquad & \text{if}\  \frac{x^2}{0.7^2}+\frac{y^2}{0.5^2}\le 1,\\
            c, \qquad &\text{otherwise},
        \end{cases}\\
        a\sim\mathcal{U}(-2, 2),\quad b\sim\mathcal{U}(0.1,2),\quad c\sim\mathcal{U}(0.5,2),\quad \omega_1 = \omega_2 = \pi.
    \end{gathered}
\end{equation}

We evaluate the proposed physics-informed framework on two-dimensional G-ROM solutions using a line-by-line reconstruction strategy, which reduces the post-processing to a sequence of one-dimensional problems.
The G-ROM baseline construction, the fine-tuning and re-training setups are detailed in Appendix~\ref{appdendix_eg3}.

\begin{figure}[!ht]
    \centering
    \subfloat[G-ROM]{\includegraphics[width=0.4\textwidth]{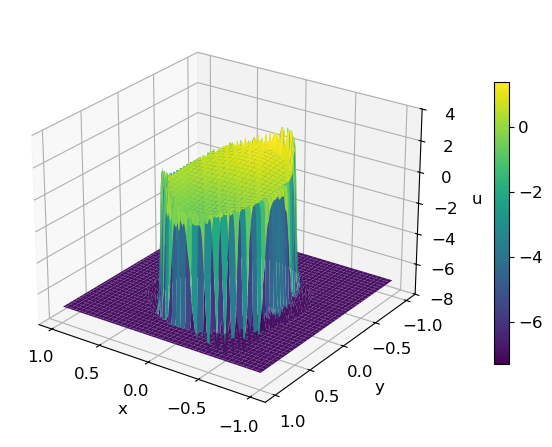}}
    \hspace{0.05\textwidth}
    \subfloat[Pre-trained model]{\includegraphics[width=0.4\textwidth]{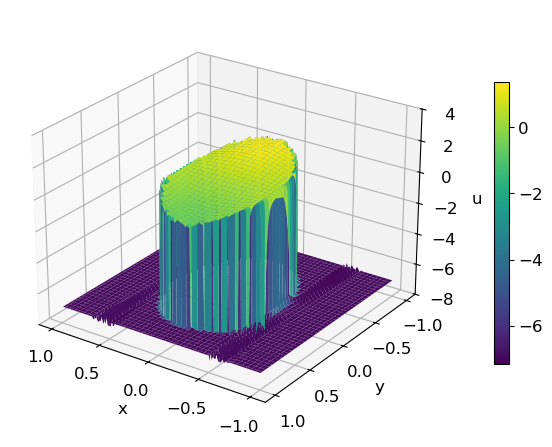}}\\
    \subfloat[Fine-tuned model]{\includegraphics[width=0.4\textwidth]{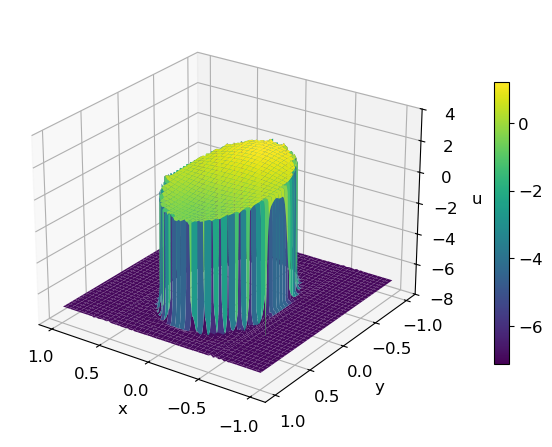}}
    \hspace{0.05\textwidth}
    \subfloat[Re-trained model]{\includegraphics[width=0.4\textwidth]{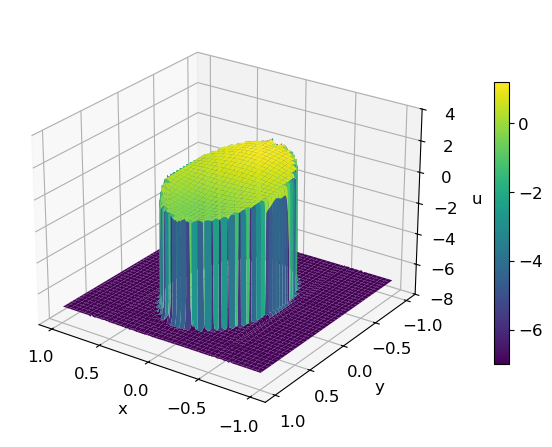}}
    \caption{Example 3. 3D view of the results before and after reconstruction, including pre-trained, fine-tuned, re-trained model.}
    \label{eg3_hyperbolic_2d_specific_sample}
\end{figure}

\begin{figure}[!ht]
    \centering
    \subfloat[Pre-trained model]{\includegraphics[width=0.3\textwidth]{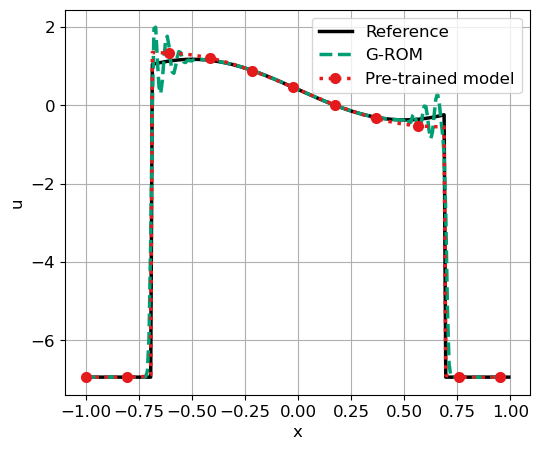}}
    \hspace{0.03\textwidth}
    \subfloat[Fine-tuned model]{\includegraphics[width=0.3\textwidth]{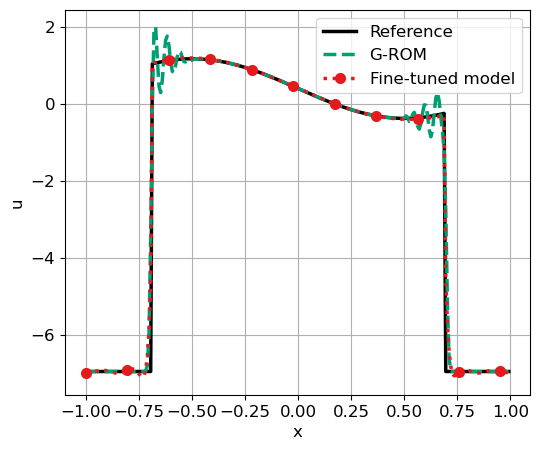}}
    \hspace{0.03\textwidth}
    \subfloat[Re-trained model]{\includegraphics[width=0.3\textwidth]{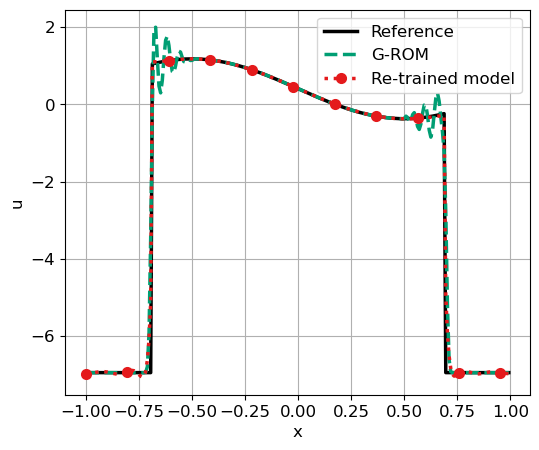}}\\
    \subfloat[Pre-trained model]{\includegraphics[width=0.3\textwidth]{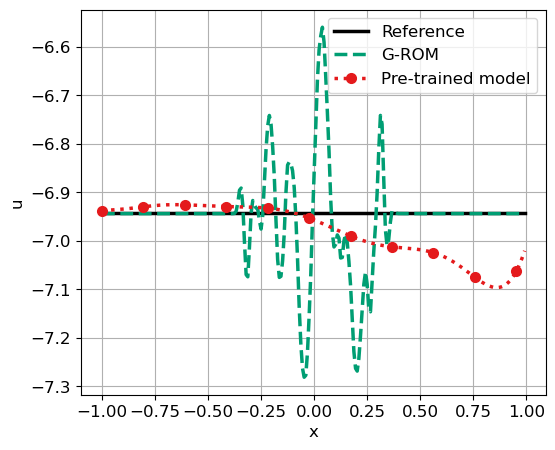}}
    \hspace{0.03\textwidth}
    \subfloat[Fine-tuned model]{\includegraphics[width=0.3\textwidth]{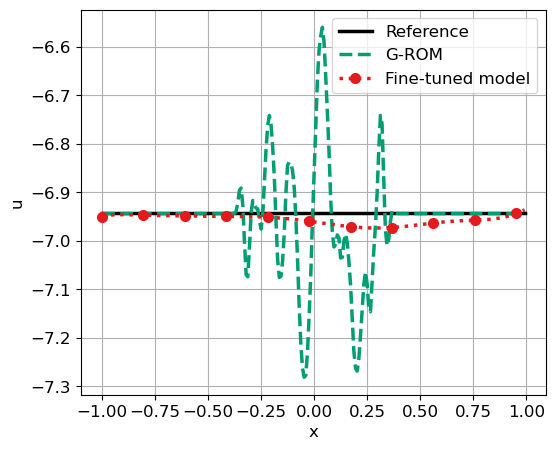}}
    \hspace{0.03\textwidth}
    \subfloat[Re-trained model]{\includegraphics[width=0.3\textwidth]{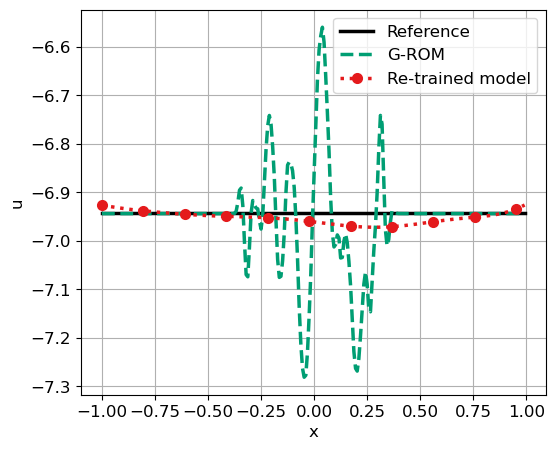}}
    \caption{Example 3. Profiles of reconstructions obtained from pre-trained, fine-tuned and re-trained model. Top: $y=-0.0625$; bottom: $y= -0.6094$ }
    \label{eg3_hyperbolic_2d_profiles}
\end{figure}

\begin{figure}
    \centering
    \subfloat[$y=-0.0625$]{\includegraphics[width=0.33\textwidth]{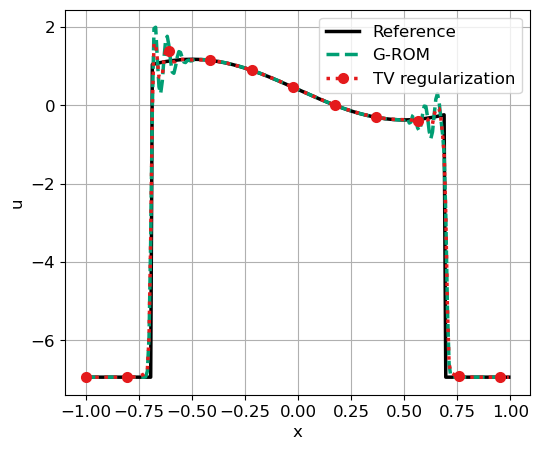}}
    \hspace{0.05\textwidth}
    \subfloat[$y=-0.6094$]{\includegraphics[width=0.33\textwidth]{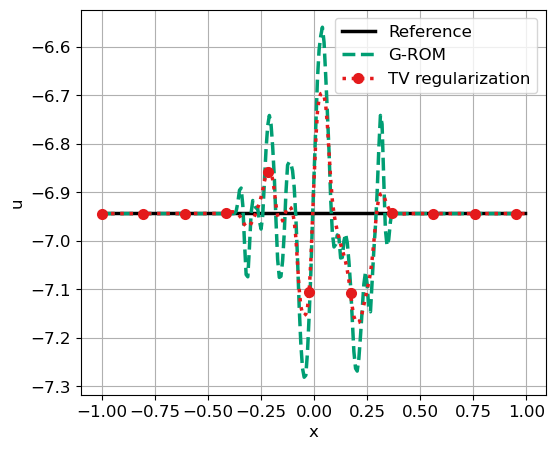}}
    \caption{Example 3. TV regularization of G-ROM solution at two selected profiles $y=-0.0625$ (left) and $y=-0.6094$ (right) with $\lambda=0.4$.}
    \label{eg3_hyperbolic_2d_specific_sample_tv}
\end{figure}

Performance is assessed on an unseen initial condition whose G-ROM solution at $T=0.5\pi$ exhibits pronounced oscillations near complex discontinuities.
Figure~\ref{eg3_hyperbolic_2d_specific_sample} shows representative reconstructions obtained with the pre-trained, fine-tuned, and re-trained models, and Table~\ref{eg3_hyperbolic_2d_specific_sample_error} reports the corresponding relative $\ell^2$ and maximum errors.
Both the fine-tuned and re-trained models substantially outperform the pre-trained model and achieve comparable accuracy. Representative one-dimensional profiles of DIOs (along $y=-0.0625$) and SOSs (along $y= -0.6094$) are shown in Figure~\ref{eg3_hyperbolic_2d_profiles}.
A comparison with total variation (TV) regularization~\cite{rudin1992nonlinear}, presented in Figure~\ref{eg3_hyperbolic_2d_specific_sample_tv}, demonstrates that while TV partially suppresses oscillations, it is less effective and underscores the advantage of the proposed framework.

\begin{table}[!ht]
    \centering
    \caption{\label{eg3_hyperbolic_2d_specific_sample_error} Example 3. Relative $\ell^2$ and maximum errors of the G-ROM solution and the reconstructions obtained for the selected test sample.}
    \begin{tabular}{c|c|c|c|c}\hline
        & G-ROM & Pre-trained model & Fine-tuned model & Re-trained model \\
        \hline
        Relative $\ell^2$& 2.6618e-02 & 1.4890e-02   &5.9820e-03  &  4.2224e-03\\
        \hline
        Maximum&1.9402e+001 & 7.1137e-01 & 4.6373e-01 & 3.7842e-01  \\
        \hline
    \end{tabular}
\end{table}

\begin{figure}[!ht]
    \centering
    \includegraphics[width=0.9\linewidth]{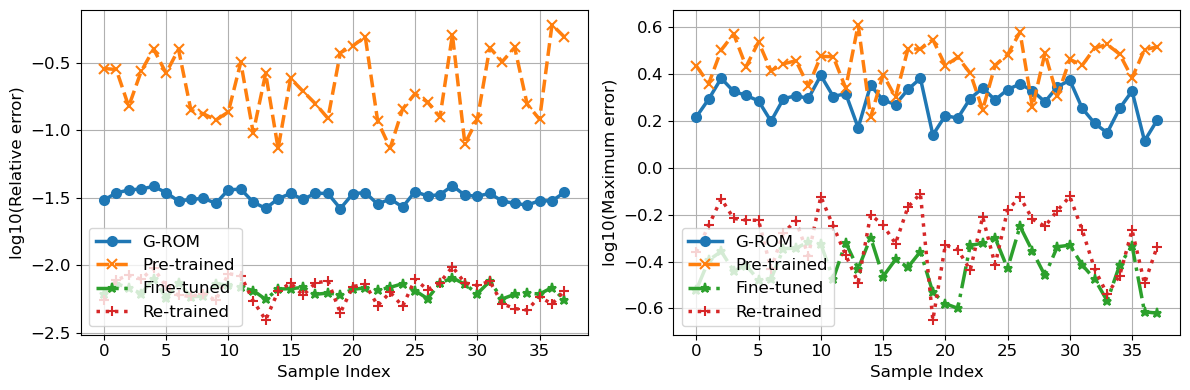}
    \caption{Example 3. Per-sample errors over $38$ G-ROM test samples for the pre-trained, fine-tuned, and re-trained models. Left: relative $\ell^2$ errors; right: maximum errors.
    }
    \label{eg3_hyperbolic_2d_model_comparision}
\end{figure}

To assess generalization, the proposed models are evaluated on $38$ newly generated test samples.
Figure~\ref{eg3_hyperbolic_2d_model_comparision} presents the distributions of the relative $\ell^2$ and maximum errors across all samples, and Table~\ref{eg3_hyperbolic_2d_model_comparision_errors} reports the corresponding statistics, including the mean, standard deviation, the $5\%$ and $95\%$ quantiles.
Both the fine-tuned and re-trained models consistently reduce errors compared with the raw G-ROM solutions, whereas the pre-trained model alone exhibits large errors and variability, underscoring the necessity of problem-specific adaptation in two-dimensional settings.
The computational cost comparison in Table~\ref{eg3_hyperbolic_2d_time_tracing} further shows that fine-tuning achieves accuracy comparable to full re-training while requiring over an order of magnitude fewer training samples and substantially reducing both dataset generation and training time.

\begin{table}[!ht]
    \centering
    \caption{\label{eg3_hyperbolic_2d_model_comparision_errors} Example 3. Statistical summary of the relative $\ell^2$ and maximum errors for the G-ROM solutions and the reconstructed solutions over $38$ test samples.
    }
    \begin{tabular}{c|l|c|c|c|c}
        \hline
        & & Mean & Std & 5\% quantile & 95\% quantile\\\hline
        \multirow{2}{*}{G-ROM} 
        & Relative $\ell^2$& 3.229e-02&3.236e-03& 2.691e-02&3.704e-02\\
        & Maximum &1.934e+00&2.984e-01&1.402e+00&2.400e+00\\\hline
        \multirow{2}{*}{Pre-trained model} 
        & Relative $\ell^2$ & 2.478e-01&1.433e-01& 7.903e-02&4.972e-01 \\
        & Maximum &2.818e+00&5.589e-01&1.797e+00&3.716e+00\\\hline
        \multirow{2}{*}{Fine-tuned model} 
        & Relative $\ell^2$ &6.579e-03&6.580e-04&5.614e-03&7.645e-03 \\
        & Maximum &3.920e-01&8.145e-02&2.503e-01&5.016e-01\\\hline
        \multirow{2}{*}{Re-trained model} 
        & Relative $\ell^2$ & 6.611e-03&1.392e-03& 4.614e-03&8.778e-03\\
        & Maximum &5.264e-01&1.426e-01&3.164e-01&7.526e-01\\\hline
    \end{tabular}   
\end{table}

\begin{table}[!ht]
    \centering
    \caption{\label{eg3_hyperbolic_2d_time_tracing}
    Example 3. Comparison of dataset generation and training costs for fine-tuned and re-trained models in the G-ROM setting for 2D hyperbolic equation. 
    }
    \begin{tabular}{c|c|c|c}\hline
         & & Dataset generation & Model training  \\ \hline
         \multirow{2}{*}{DIOs} 
         &Fine-tuned model & 85.43 s / 230 samples & 11.82 s / 200 epochs\\ 
         &Re-trained model & 836.11 s / 2152 samples & 317.30 s / 600 epochs\\\hline
         \multirow{2}{*}{SOSs} 
         &Fine-tuned model & 65.78 s / 232 samples & 10.83 s / 200 epochs\\ 
         &Re-trained model & 655.12 s / 2320 samples & 375.61 s / 600 epochs\\\hline
    \end{tabular}
\end{table}

\subsubsection{Example 4: 2D inviscid Burgers' equation \label{eg4_burgers_2d}}

Consider the 2D Burgers' equation
\begin{equation}
    \label{eq_burgers_2d}
    u_t + uu_x+uu_y = 0 , \qquad (x,y)\in [0,4]\times[0,4],
\end{equation}
with periodic boundary condition in each direction and the parametric initial condition 
\begin{equation}
    \label{eq_burgers_2d_ic}
    \begin{gathered}
        u(x,y,0) = u_0(x,y)= a+b\sin(\omega(x+y)+\varphi),\\
        \qquad a\sim\mathcal{U}(-2,2),\quad  
        b\sim\mathcal{U}(0.1,2),\quad  
        \varphi\sim\mathcal{U}(-\pi,\pi),\quad
        \omega=\pi.
    \end{gathered}
\end{equation}

For this two-dimensional nonlinear example, the G-ROM solution exhibits oscillatory artifacts near moving discontinuity interfaces. 
We apply the proposed reconstruction in a line-by-line manner, with both fine-tuning and re-training datasets generated from G-ROM solutions. Details of the G-ROM baseline construction, the fine-tuning and re-training setups are provided in Appendix~\ref{appendix_eg4}.

\begin{figure}[!ht]
    \centering
    \subfloat[Pre-trained model]{\includegraphics[width=0.3\textwidth]{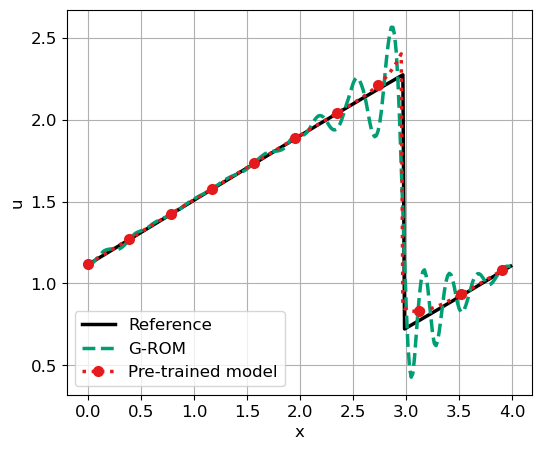}}
    \hspace{0.03\textwidth}
    \subfloat[Fine-tuned model]{\includegraphics[width=0.3\textwidth]{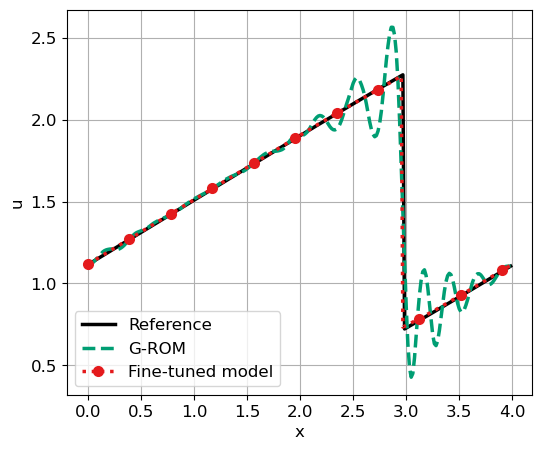}}
    \hspace{0.03\textwidth}
    \subfloat[Re-trained model]{\includegraphics[width=0.3\textwidth]{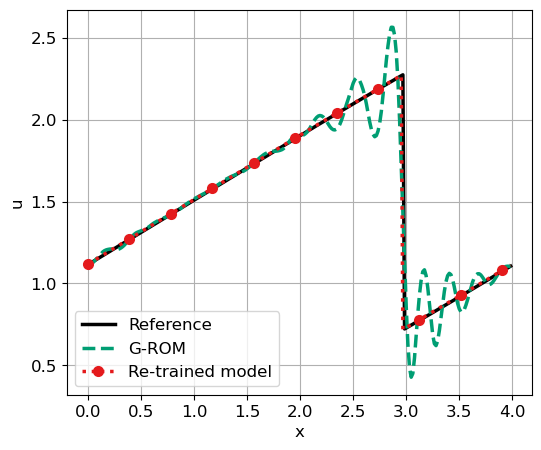}}
    \caption{Example 4. Profiles at $y=2.00$ of reconstructions obtained from pre-trained, fine-tuned, and re-trained model.}
    \label{eg4_burgers_2d_profiles}
\end{figure}

Figure~\ref{eg4_burgers_2d_specific_sample} and Table~\ref{eg4_burgers_2d_specific_sample_error} compare the reconstructions obtained by the pre-trained, fine-tuned, and re-trained models for an unseen test sample at $T=1.0$. 
A representative cross-section at $y=0.5$ is shown in Figure~\ref{eg4_burgers_2d_profiles}. 
The pre-trained model provides only limited correction, whereas the fine-tuned and re-trained models substantially suppress oscillations and recover sharper solution profiles.

\begin{figure}[!ht]
    \centering
    \subfloat[G-ROM]{\includegraphics[width=0.4\textwidth]{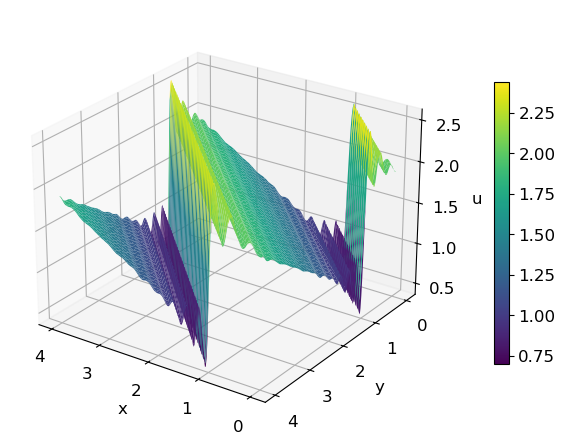}}
    \hspace{0.05\textwidth}
    \subfloat[Pre-trained model]{\includegraphics[width=0.4\textwidth]{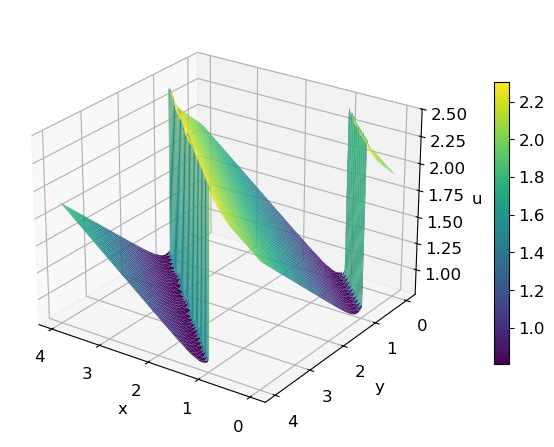}}\\
    \subfloat[Fine-tuned model]{\includegraphics[width=0.4\textwidth]{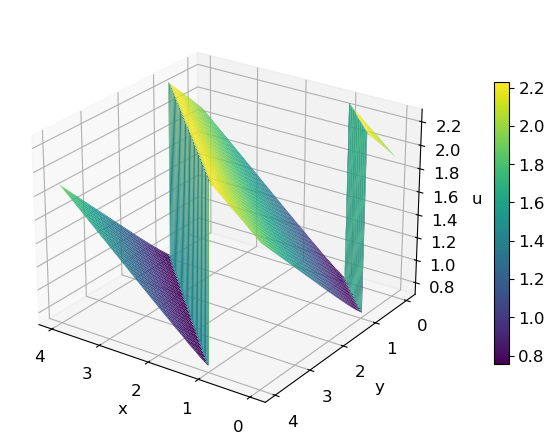}}
    \hspace{0.05\textwidth}
    \subfloat[Re-trained model]{\includegraphics[width=0.4\textwidth]{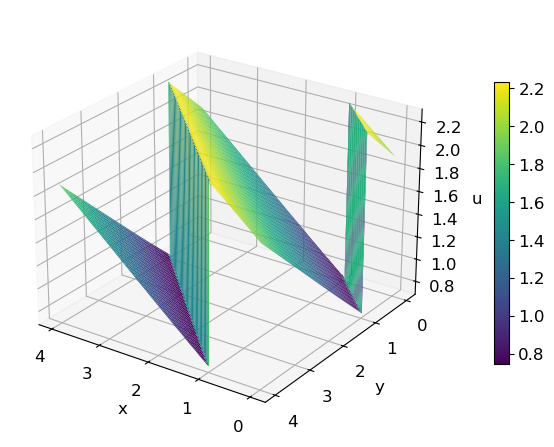}}
    \caption{Example 4. 3D view of the results before and after reconstruction, including pre-trained, fine-tuned, and re-trained model.}
    \label{eg4_burgers_2d_specific_sample}
\end{figure}

\begin{table}[!ht]
    \centering
    \caption{\label{eg4_burgers_2d_specific_sample_error} Example 4. 
    Relative $\ell^2$ and maximum errors of the G-ROM solutions and the reconstructions for the selected test sample.}
    \begin{tabular}{c|c|c|c|c}\hline
        & G-ROM & Pre-trained model & Fine-tuned model & Re-trained model \\
        \hline
        Relative $\ell^2$& 5.9726e-02 & 1.4677e-02   &2.7290e-03   &  1.9985e-03\\
        \hline
        Maximum&3.5203e-01 & 1.5876e-01 & 2.0617e-02 & 2.3885e-02 \\
        \hline
    \end{tabular}
\end{table}

\begin{figure}[!ht]
    \centering
    \includegraphics[width=0.9\linewidth]{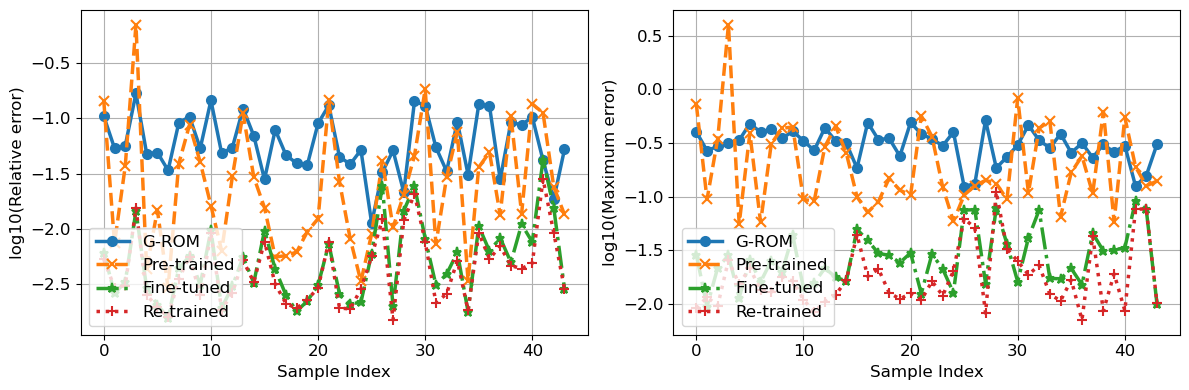}
    \caption{Example 4. Per-sample errors over $44$ G-ROM test samples for the pre-trained, fine-tuned, and re-trained models. Left: relative $\ell^2$ errors; right: maximum errors.}
    \label{eg4_burgers_2d_model_comparision}
\end{figure}

To assess robustness, we further evaluate the models on $44$ newly generated test samples with intersecting characteristics. 
Table~\ref{eg4_burgers_2d_model_comparision_errors} shows that fine-tuning reduces both the mean relative and maximum error by approximately one order of magnitude, with per-sample error curves shown in Figure~\ref{eg4_burgers_2d_model_comparision}.
The performance of fine-tuning is close to that of the fully re-trained model, confirming that problem-specific adaptation remains effective for two-dimensional nonlinear transport dynamics.

\begin{table}[!ht]
    \centering
    \caption{\label{eg4_burgers_2d_model_comparision_errors} Example 4. 
    Statistical summary of the relative $\ell^2$ and maximum errors for the G-ROM solutions and the reconstructed solutions over $44$ test samples.
    }
    \begin{tabular}{c|l|c|c|c|c}
        \hline
        & & Mean & Std & 5\% quantile & 95\% quantile\\\hline
        \multirow{2}{*}{G-ROM} 
        & Relative $\ell^2$& 6.995e-02&4.002e-02& 2.237e-02&1.418e-01\\
        & Maximum &3.203e-01&9.741e-02&1.361e-01&4.877e-01\\\hline
        \multirow{2}{*}{Pre-trained model} 
        & Relative $\ell^2$ & 5.550e-02&1.093e-01&3.500e-03&1.456e-01 \\
        & Maximum &3.300e-01&5.926e-01&3.300e-01&5.926e-01\\\hline
        \multirow{2}{*}{Fine-tuned model} 
        & Relative $\ell^2$ &7.045e-03&7.578e-03&1.799e-03&2.298e-02 \\
        & Maximum &3.172e-02&2.100e-02&1.162e-02&7.616e-02\\\hline
        \multirow{2}{*}{Re-trained model} 
        & Relative $\ell^2$ & 5.608e-03&5.239e-03& 1.814e-03&1.493e-02\\
        & Maximum &2.301e-02&2.147e-02&8.452e-03&7.436e-02\\\hline
    \end{tabular}   
\end{table}

It is worth emphasizing that the dominant online cost comes from evaluating Gegenbauer reconstructions over candidate parameter pairs along each spatial line for this two-dimensional problem.
As a result, the per-sample reconstruction time is nearly identical for the pre-trained, fine-tuned, and re-trained models, with runtimes concentrated around $47.5\,\mathrm{s}$. 
Thus, the choice of parameter predictor mainly affects reconstruction accuracy rather than online evaluation cost.

From the offline computational perspective, Table~\ref{eg4_burgers_2d_time_tracing} shows that fine-tuning is substantially more efficient than full re-training. 
Although the re-trained model can be slightly more accurate in some cases, fine-tuning achieves comparable reconstruction accuracy with more than an order-of-magnitude reduction in both dataset generation and model training costs. 
These results demonstrate that the proposed two-stage strategy provides an efficient and practical approach for high-quality oscillation correction in two-dimensional settings.

\begin{table}[!ht]
    \centering
    \caption{\label{eg4_burgers_2d_time_tracing} Example 4. 
    Comparison of dataset generation and training costs for fine-tuned and re-trained models in the G-ROM setting for 2D Burgers' equation.}
    \begin{tabular}{c|c|c}\hline
         & Dataset generation & Model training  \\ \hline
         Fine-tuned model & 77.57 s / 150 samples & 11.37 s / 200 epochs\\ \hline
         Re-trained model & 1237.19 s / 2204 samples & 340.32 s / 600 epochs\\\hline
    \end{tabular}
\end{table}

\subsection{\texorpdfstring{{Fine-tuning Analysis: Accuracy, Efficiency, and Robustness}}{Fine-tuning Analysis: Accuracy, Efficiency, and Robustness}\label{Accuracy_and_efficiency_tradeoff}}

In practical applications of physics-informed reconstruction, there is an inherent trade-off between reconstruction accuracy and computational cost. 
While full re-training often provides a strong problem-specific baseline, lighter-weight approaches, such as fine-tuning, offer a more efficient alternative whose performance depends on the available target data, the choice of trainable components, and the robustness of the pre-trained representation. 
This dependence is consistent with broader observations in modern deep learning, where generalization performance can vary nontrivially with data availability, model capacity, and training strategy~\cite{nakkiran2021deep}.

Using the one-dimensional Burgers' equation as a representative nonlinear benchmark, we conduct a problem-driven analysis of the proposed two-stage strategy. 
Specifically, we examine three key aspects: the amount of data required for effective fine-tuning, the choice of network layers to be fine-tuned, and the robustness of the adapted model under distributional shifts.
These analyses clarify the accuracy-efficiency trade-off and guide the practical use of the proposed probabilistic Gegenbauer reconstruction framework.

\subsubsection{Accuracy and efficiency versus fine-tuning dataset volume\label{Accuracy and computational cost versus fine-tuning dataset volume}}

The effectiveness of fine-tuning depends strongly on the amount of problem-specific training data. 
Insufficient data may lead to overfitting and suboptimal adaptation, whereas large samples generally improve reconstruction accuracy at the expense of additional data-generation and training cost. 
Such data-dependent behavior is also commonly observed in supervised domain adaptation, where additional in-domain samples often lead to improved but diminishing gains~\cite{miceli2017regularization}. 
We therefore examine how the fine-tuning dataset size affects the accuracy-efficiency trade-off of the proposed parameter-selection model.

Fine-tuning datasets are constructed from G-ROM solutions of the one-dimensional Burgers' equation using the setup in Subsection~\ref{eg2_1d_burgers}. 
We consider
\begin{equation*}
    N_{\mathrm{data}} = [5, 8, 18, 39, 76, 161, 247].
\end{equation*}
All models are evaluated on the same fixed Burgers G-ROM test set.

The computational cost increases with the fine-tuning dataset size, with dataset generation forming the dominant contribution because each additional sample requires a new G-ROM solution. 
For small datasets, the training and validation losses show a noticeable gap, indicating overfitting; this gap decreases as more target-specific samples are included. 
Detailed cost and loss diagnostics are provided in Appendix~\ref{appendix_finetuning_data_volume}.

Figure~\ref{statistic_loss_vs_volume} reports the mean errors and empirical $5\%$-$95\%$ quantile ranges over the Burgers' G-ROM test set as the fine-tuning dataset size increases, with the corresponding per-sample reconstruction error curves complemented in Figure~\ref{accuracy_and_cost_vs_volume}. 
The errors decrease rapidly when a small number of target-specific samples are added and then gradually approach the accuracy of the fully re-trained model. 
The observed decay is close to the Monte Carlo reference rate $O(N_{\mathrm{data}}^{-1/2})$~\cite{caflisch1998monte}, suggesting that the adapted model has sufficient capacity over the tested range, and the remaining reconstruction error is mainly governed by sampling variability rather than optimization failure or severe model underfitting.

Indeed, the fine-tuning task amounts to learning the input-dependent map
\begin{equation*}
    v(x) \mapsto p_\theta(m,\lambda\mid v),
\end{equation*}
where different shock locations, oscillation amplitudes, and local geometries may favor different Gegenbauer parameter distributions. 
When $N_{\mathrm{data}}$ is small, the target-specific samples provide limited coverage of these local configurations, leading to biased or unstable parameter predictions. 
As $N_{\mathrm{data}}$ increases, the coverage improves and the prediction error decreases, but the gains exhibit diminishing returns consistent with sampling-error reduction.

\begin{figure}[!ht]
    \centering
    \includegraphics[width=0.75\linewidth]{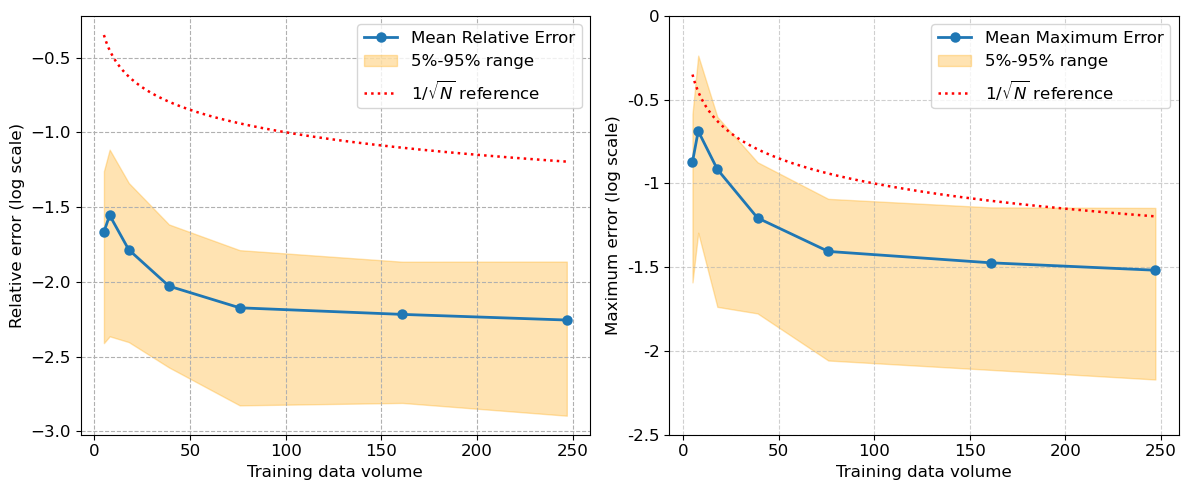}
    \caption{Statistical analysis of reconstruction errors over 38 test samples as a function of fine-tuning dataset size. Shaded regions indicate the 5\%–95\% quantile range, and the mean error is plotted for both relative $\ell^2$ (left) and maximum error (right). A $1/\sqrt{N}$ reference is included for comparison.}
    \label{statistic_loss_vs_volume}
\end{figure}

\subsubsection{Accuracy and efficiency versus fine-tuning strategies\label{network_expressivity}}

The effectiveness and efficiency of fine-tuning depend not only on the amount of problem-specific data, but also on which components of the pre-trained network are updated. 
To examine this dependence, we conduct an ablation study over several adaptation strategies, ranging from updating only the final fully-connected layer head to fine-tuning all convolutional and dense modules.
The detailed backward-pass cost estimates, which quantify the computational overhead of updating different layer groups, are reported in Appendix Table~\ref{tab_finetune_grad_formula}.

Table~\ref{tab_network_expressivity} compares different fine-tuning strategies in terms of trainable-parameter ratio, training time, and optimal training / validation losses. 
Fine-tuning only the fully-connected prediction head achieves a substantial loss reduction with minimal overhead. 
Updating the fusion convolution provides additional improvement but increases the number of trainable parameters by approximately a factor of $3.5$. 
Full fine-tuning gives the lowest training loss, yet the improvement over lighter strategies is modest compared with its additional cost. 
In contrast, strategies that fine-tune only parallel convolutional or fusion layers alone perform noticeably worse.

\begin{table}[!ht]
    \centering
    \footnotesize
    \caption{Ablation study of fine-tuning strategies for \texttt{CNN\_ReconParamNet}, reporting the proportion of trainable parameters, training time, and training / validation losses.}
    \label{tab_network_expressivity}
    \begin{tabular}{c|c|c|c|c}
        \hline
        \textbf{Strategy}& \textbf{Trainable parameters} & \textbf{Training time} & \textbf{Training loss} & \textbf{Validation loss} \\\hline
        \textbf{only fc}& $42,460/216,412\approx19.62\%$ & 9.81\ s & 3.0460e-03 & 3.4249e-03\\ \hline
        \textbf{fc+fusion} & $153,628/216,412\approx 70.99\%$ & 10.47\ s & 2.4077e-03& 3.2686e-03\\ \hline
        \textbf{fc+fusion+large} & $182,876/216,412\approx 84.50\%$ & 11.39\ s & 2.4272e-03& 3.3637e-03\\ \hline
        \textbf{only fusion}& $111,168/216,412\approx 51.37\%$ & 10.27\ s & 6.1843e-03 & 6.6092e-03\\ \hline
        \textbf{only conv}& $62,784/216,412\approx 29.01\%$ & 11.23\ s & 3.5949e-02 & 3.0451e-02\\ \hline
        \textbf{all} & $216,412/216,412= 100.00\%$ & 13.51\ s & 2.3858e-03 & 3.8405 e-03\\ \hline
    \end{tabular}
\end{table}

\begin{figure}[!ht]
    \centering
    \includegraphics[width=0.98\linewidth]{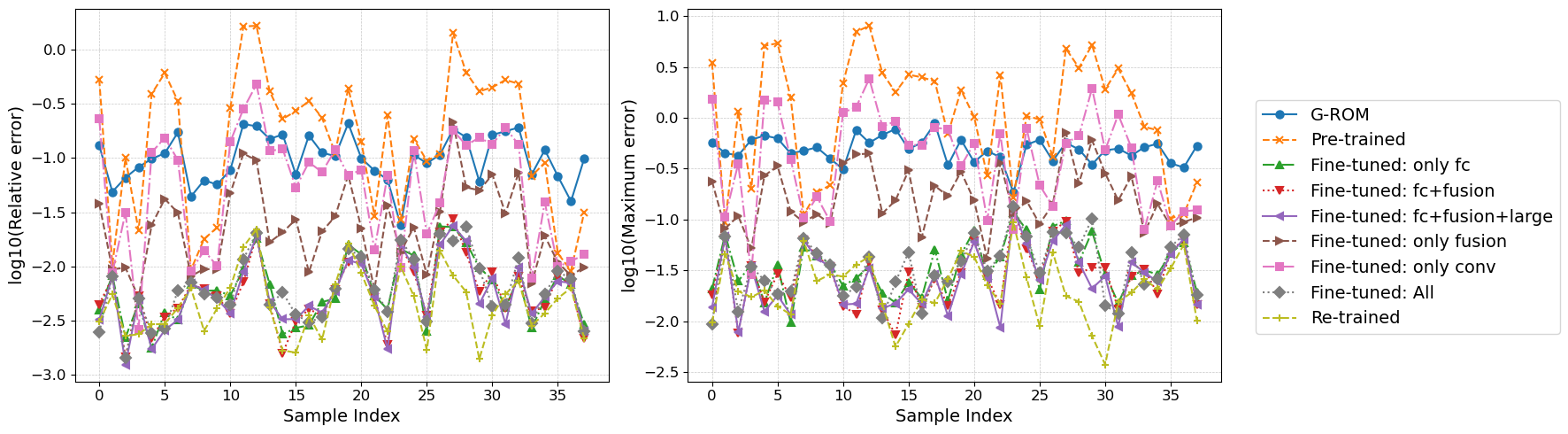}
    \caption{Reconstruction accuracy comparison on $38$ Burgers' G-ROM test samples. Left: relative $\ell^2$ error; right: maximum pointwise error. 
    Results are reported for the pre-trained model, fine-tuned models under different fine-tuning strategies, and the fully retrained model.}
    \label{fig_network_expressivity}
\end{figure}

Figure~\ref{fig_network_expressivity} reports the reconstruction errors over the $38$ Burgers' G-ROM test samples for different fine-tuning strategies. 
Fine-tuning only the final fully-connected prediction head, which accounts for less than $20\%$ of the network parameters, already achieves accuracy comparable to more aggressive and extensive fine-tuning strategies. 
Expanding the trainable subset yields only marginal additional improvement, whereas updating only convolutional feature-extraction modules leads to inferior performance. 
These results indicate that the convolutional backbone learned during pre-training provides transferable oscillation-aware features, and that task-specific adaptation primarily requires recalibrating the high-level map from extracted features to Gegenbauer parameter distributions.

\subsubsection{Out-of-distribution generalization}

Out-of-distribution (OoD) robustness is important for assessing whether a fine-tuned reconstruction model remains reliable beyond the data used for adaptation. 
Here we investigate this behavior with two representative distributional shifts: temporal shifts and additive perturbations of the G-ROM input. 
This complements existing analyses of OoD behavior in linear operator learning~\cite{de2023convergence} and provides a nonlinear test of the proposed physics-informed reconstruction framework.

For temporal OoD tests, we follow the Burgers G-ROM setup but evaluate the models at randomly sampled times $t\in[0.6,1.2]$. 
For noise OoD tests, additive Gaussian perturbations are applied to the G-ROM input, 
\begin{equation} 
    \sigma_{\mathrm{noise}} \sim \mathcal{U}(0.01,0.1)\bigl(\max(u)-\min(u)\bigr), \qquad u_T^{\mathrm{noisy}} = u_T+\mathcal{N}(0,\sigma_{\mathrm{noise}}^2), 
\end{equation} 
where $u_T$ denotes the G-ROM solution at the test time. 
Detailed per-sample OoD reconstruction results are provided in Appendix~\ref{appedix_ood}.

\begin{figure}[!ht]
    \centering
    \includegraphics[width=0.75\linewidth]{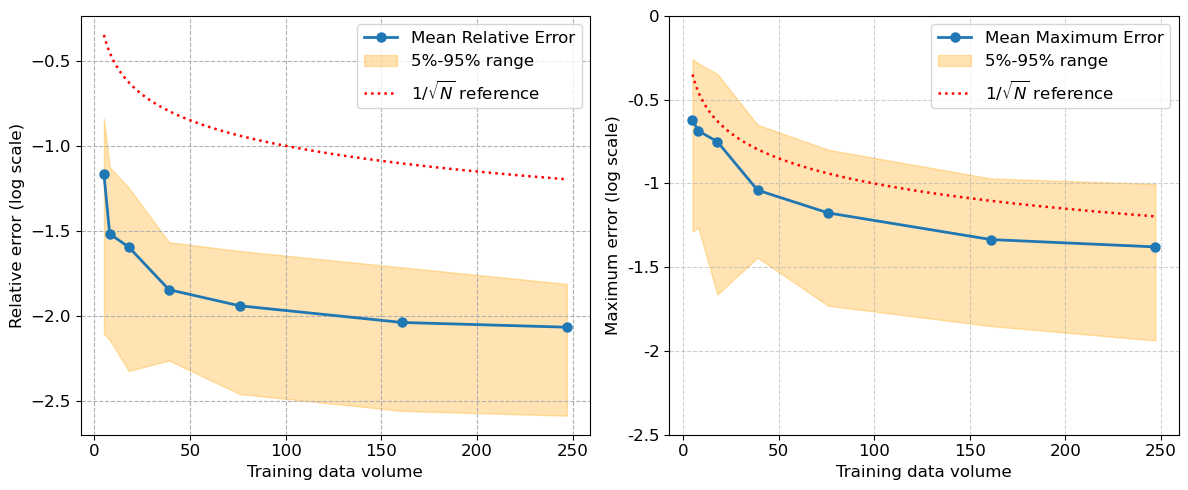} 
    \caption{Statistical analysis of reconstruction errors under time-based out-of-distribution perturbations as a function of fine-tuning dataset volume. Shaded regions indicate 5\%–95\% quantiles; mean errors are shown for relative $\ell^2$ (left) and maximum errors (right).}
    \label{time_OoD_data_volume_loss}
\end{figure}

\begin{figure}[!ht]
    \centering
    \includegraphics[width=0.75\linewidth]{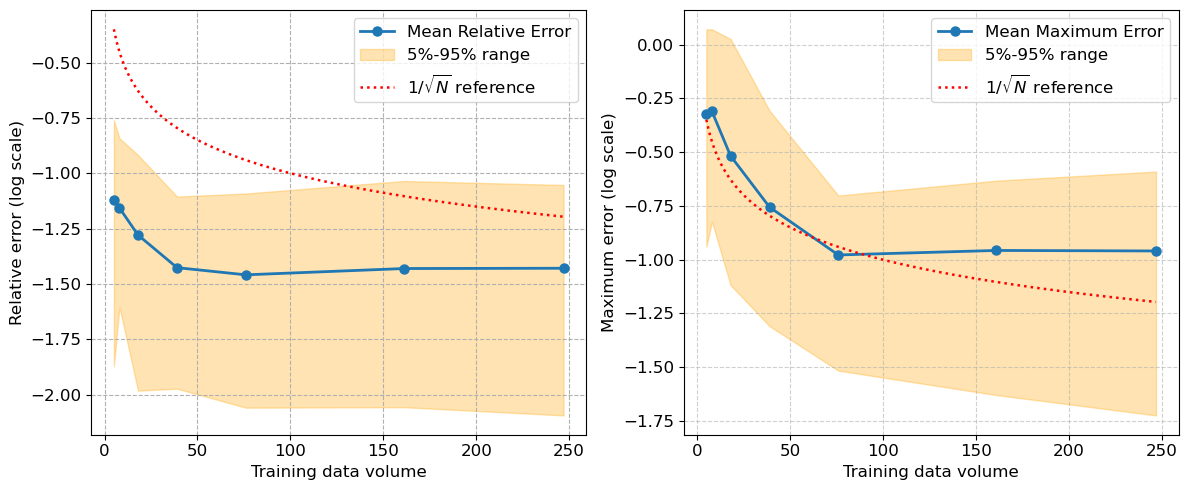} 
    \caption{Statistical analysis of reconstruction errors under noise-based out-of-distribution perturbations as a function of fine-tuning dataset volume. Shaded regions indicate 5\%–95\% quantiles; mean errors are shown for relative $\ell^2$ (left) and maximum errors (right).}
    \label{nosie_OoD_data_volume_loss}
\end{figure}

To quantify how target-specific data volume affects OoD generalization, we evaluate models fine-tuned with different dataset sizes on $41$ temporal-OoD samples and $36$ noise-OoD samples, and Figures~\ref{time_OoD_data_volume_loss} and~\ref{nosie_OoD_data_volume_loss} report the mean errors and empirical $5\%$-$95\%$ quantile ranges as functions of the fine-tuning dataset size. For temporal OoD tests, the errors decrease consistently as the fine-tuning dataset size increases and approximately follow the Monte Carlo reference rate $O(N_{\mathrm{data}}^{-1/2})$, suggesting that temporal generalization is mainly limited by the sampling coverage of target-specific configurations, such as shock locations and oscillation patterns. 
In contrast, for noise OoD tests, the error reduction saturates beyond a moderate dataset size and remains above the Monte Carlo reference trend. This indicates that robustness to input perturbations is not controlled solely by the amount of fine-tuning data and may require additional denoising or robustness-aware training mechanisms.

    \section{Conclusions and Future Work \label{conclusion}}

This paper proposed a physics-informed machine-learning framework for post-processing oscillatory data-driven solutions of transport-dominated problems.
The framework predicts a probability distribution over candidate Gegenbauer reconstruction parameters and uses it to form a weighted reconstruction that accounts for parameter uncertainty.
The method was validated on one- and two-dimensional transport-dominated problems using G-ROM and DeepONet as representative data-driven solvers, and the reconstruction reduced errors by up to one to two orders of magnitude in the tested cases. 
By fine-tuning a pre-trained predictor with a small problem-specific dataset, the two-stage strategy achieved accuracy comparable to full re-training at substantially lower computational cost.
The fine-tuning analysis showed that modest target-specific datasets already provide substantial gains, and that updating only the final prediction head offers an effective accuracy-efficiency compromise.
The out-of-distribution tests further showed that fine-tuning improves robustness under temporal shifts and noisy inputs, although noise perturbations remain more challenging.
Overall, these results demonstrate that the proposed framework provides an accurate, efficient, and flexible post-processing tool for data-driven solvers of transport-dominated systems.

Several directions remain for future work.
A natural extension is to apply the framework to more complex nonlinear and multiscale PDEs across broader physical regimes.
Incorporating uncertainty quantification for reconstruction errors could provide confidence estimates and further improve reliability.
Adaptive or active fine-tuning strategies that identify informative target samples could further reduce data-generation costs.
Finally, transfer learning across different classes of PDEs may enable more efficient reuse of learned representations and accelerate adaptation to new problem settings.

    \appendix

\section*{Appendix}
\addcontentsline{toc}{section}{Appendix}

\section{Details for Dataset Generation during Pre-training \label{appendix_dataset_pretraining}}

All samples are generated on the same uniform grid as the discontinuity-induced dataset, with $n_{\mathrm{space}}=256$ over $x\in[-1,1)$.
And for each profile, Gegenbauer reconstructions are evaluated over the candidate set spanned with $m=0,\ldots,10$ and $\lambda=1,\ldots,20$. 
The pre-training data consist of two classes of oscillatory profiles: discontinuity-induced oscillations and spurious numerical oscillations.

\subsection{Dataset with Discontinuity-Induced Oscillations}
\label{appendix_dataset_true}

\paragraph{Problem setup:}
Training samples are generated from a representative one-dimensional linear transport problem (Subsection~\ref{example1}) within a reduced-order modeling framework to capture physically meaningful oscillations near discontinuities.

\paragraph{Projection-based low-cost construction:}

To avoid repeatedly solving the full Galerkin reduced system during large-scale data generation, we use a projection-based construction. 
A snapshot matrix is assembled from randomly shifted smooth polynomial profiles and compressed by truncated singular value decomposition. 
The leading left singular vectors define a reduced space, onto which additional shifted polynomial samples are projected to produce low-cost oscillatory profiles.

\paragraph{Parameter settings:} 
The snapshot matrix is constructed from $500$ randomly sampled polynomial initial conditions of degree at most $10$, evolved over $50$ time levels on $t\in[0,1]$. Before compression, the snapshots are smoothed using a Gaussian filter with width $\sigma=2.0$. The reduced space used for data generation has dimension $r_{\mathrm{proj}}=15$. Using this basis, $2000$ projected solutions are generated at the target time $T=0.2$. The total dataset-generation time for this class is $370.42\,\mathrm{s}$.

\subsection{Dataset with Spurious Oscillations}
\label{appendix_dataset_spurious}

\paragraph{Baseline profile:} 
The baseline field $u_{\mathrm{base}}(x)$ is chosen as a constant or a low-degree polynomial with degree $d\le d_{\max}$.

\paragraph{Localized oscillatory perturbation:}

A localized oscillation is generated by a sinusoid modulated by a Gaussian envelope,
\begin{equation}
    \mathcal{O}(x) = A \exp\!\left(-\frac{(x-x_c)^2}{2\sigma_w^2}\right) \sin(\omega x),
\end{equation}
where $x_c$ and $\sigma_w$ denote the center and width of the envelope, respectively. 
The amplitude $A$ and the dimensionless frequency parameter $f$ are sampled independently, e.g. $A\sim\mathcal{U}(A_{\min},A_{\max})$ and $f\sim\mathcal{U}(f_{\min},f_{\max})$.

\paragraph{Smooth stochastic perturbation:}
A localized smooth-noise component is constructed by filtering a white-noise field $\eta(x)$ with a Gaussian kernel $G_{\sigma_n}$ and applying the same localization envelope:
\begin{equation}
    \mathcal{N}(x) = \varepsilon \exp\!\left(-\frac{(x-x_c)^2}{2\sigma_w^2}\right) \widetilde{\eta}(x),
    \qquad
    \widetilde{\eta}(x) = G_{\sigma_n}\ast \eta(x),
\end{equation}
where $\varepsilon\sim\mathcal{U}(\varepsilon_{\min},\varepsilon_{\max})$ controls the noise amplitude.

\paragraph{Parameter settings:}

In the numerical experiments, the baseline is chosen as a constant field. 
The oscillation amplitude, dimensionless frequency, envelope width, and noise amplitude are sampled as
\begin{equation*}
    A\sim\mathcal{U}(0.1,0.4),\qquad
    f\sim\mathcal{U}(2,7),\qquad
    \sigma_w\sim\mathcal{U}(0.05,0.2),\qquad
    \varepsilon\sim\mathcal{U}(0.01,0.1),
\end{equation*}
with $x_c=0$ and $\omega=\pi f$. 
A total of $2000$ spurious samples are generated.

\section{Detailed Pre-training Network Architecture \label{pretrain_network}}
    
Table~\ref{tab_pretrain_network} shows the complete layer-wise description of the \texttt{CNN\_ReconParamNet} network, including output shapes and parameter counts.

\begin{table}[!ht]
    \centering
    \small
    \caption{Summary of the \texttt{CNN\_ReconParamNet} architecture.}
    \label{tab_pretrain_network}
    \begin{tabular}{c|l|c|l}
    \hline
    \textbf{Module} & \textbf{Main operation} & \textbf{Output shape} & \textbf{Parameters} \\
    \hline
    \texttt{conv\_small}  
    & Two Conv1d layers, kernel size $3$ 
    & $[1,64,256]$ 
    & $12{,}608$ \\
    \hline
    \texttt{conv\_medium} 
    & Two Conv1d layers, kernel size $5$ 
    & $[1,64,256]$ 
    & $20{,}928$ \\
    \hline
    \texttt{conv\_large}  
    & Two Conv1d layers, kernel size $7$ 
    & $[1,64,256]$ 
    & $29{,}248$ \\
    \hline
    \texttt{fusion\_conv} 
    & Conv1d fusion + BatchNorm1d 
    & $[1,192,256]$ 
    & $111{,}168$ \\
    \hline
    Pooling 
    & Adaptive average pooling 
    & $[1,192]$ 
    & $0$ \\
    \hline
    \texttt{fc} 
    & Linear$(192,220)$ 
    & $[1,220]$ 
    & $42{,}460$ \\
    \hline
    \multicolumn{3}{r|}{\textbf{Total trainable parameters}} 
    & $216{,}412$ \\
    \hline
    \multicolumn{3}{r|}{\textbf{Total mult-adds}} 
    & $44.48$ MB \\
    \hline
    \end{tabular}
\end{table}

\section{Additional Details for Linear Transport Problem\label{appendix_eg1}}

\subsection{G-ROM baseline construction\label{appendix_eg1_grom}}

$100$ randomly sampled polynomial initial conditions of degree at most $10$ on a uniform grid with $n_x=256$ points are generated.
For each initial condition, snapshots are collected at $n_t=50$ time instances using the exact shifted solution. 
The snapshot data are smoothed with a Gaussian filter of width $\sigma=2.0$ and compressed by singular value decomposition, with leading $r=20$ POD modes retained to form the reduced space.

\subsection{Problem-specific adaptation setup\label{appendix_eg1_adaptation_setup}}

The fine-tuning dataset contains $100$ valid G-ROM samples at $T=0.2$, while the re-training dataset contains $2000$ samples. 
The fine-tuned model is initialized from the corresponding pre-trained predictor, with the convolutional backbone frozen and only the final prediction head updated for $200$ epochs using batch size $16$ and a cosine-annealing learning rate from $10^{-4}$ to $10^{-5}$. 
The re-trained model uses the same architecture but is trained from scratch on the larger problem-specific dataset.

\section{Additional Details for Inviscid Burgers' Problem\label{appendix_eg2}}

\subsection{G-ROM setups and additional visualization results}
\label{appendix_eg2_grom}

\paragraph{G-ROM baseline construction:}

The reduced-order space is constructed from snapshots generated by evolving $100$ random initial conditions at $50$ uniformly distributed time instances over $t\in[0,1]$ using the method of characteristics. 
To improve the decay of the Kolmogorov $n$-width, all snapshots are preprocessed with a Gaussian filter of width $\sigma=1$, and the leading $r=30$ POD modes obtained from SVD are retained to define the G-ROM space.

\paragraph{Problem-specific adaptation setup:}

The fine-tuning dataset contains $158$ valid samples generated from G-ROM solutions at $T=0.6$, with each sample evolved from a randomly prescribed initial condition using a third-order Runge-Kutta scheme with $500$ time steps, while the re-training dataset contains $1926$ samples.
During fine-tuning, the convolutional backbone is frozen and only the fully-connected layer is updated for $200$ epochs with batch size $16$ and a cosine-annealing learning-rate schedule from $10^{-3}$ to $10^{-5}$. 
For comparison, the re-trained model uses the same architecture and loss formulation but is trained from scratch on the larger dataset, with all network parameters updated.

\subsection{DeepONet setups and additional visualization results\label{appendix_eg2_deeponet}}

\paragraph{DeepONet baseline construction:}

The DeepONet approximates the solution operator from the initial condition of the inviscid Burgers' equation to the solution at $T=0.6$. 
The branch network takes the discretized initial condition, augmented with its parameters, as input, while the trunk network encodes spatial coordinates using periodic Fourier features~\eqref{eq:deeponet_periodic_encoding} with $K=3$. 
Both networks are fully-connected with three hidden layers of width $256$ and $\tanh$ activations. 
The model is trained with Adam, an initial learning rate of $10^{-3}$, and a cosine-annealing schedule for $2\times10^4$ epochs, minimizing the relative $\ell^2$ error.

\paragraph{Problem-specific adaptation setup:}

Fine-tuning and re-training datasets are generated from DeepONet predictions at $T=0.6$ for randomly sampled initial conditions, retaining only samples with shock formation, yielding $134$ fine-tuning samples and $1904$ re-training samples. 
The parameter predictor is then fine-tuned or re-trained using the same architecture and training strategy as in the G-ROM setting. 
During fine-tuning, the convolutional backbone is frozen and only the final prediction head is updated, whereas the re-trained model is trained from scratch on the larger DeepONet-specific dataset.

\section{Additional Details for 2D Linear Equation \label{appdendix_eg3}}

\subsection{G-ROM baseline construction}

The spatial domain is discretized on a $256\times256$ uniform grid.
A reduced-order space of dimension $r=30$ is constructed by applying SVD to a snapshot matrix generated from $50$ random initial conditions, each evolved from $t=0$ to $2\pi$ at $100$ time instances using exact solutions and preprocessed with Gaussian filtering ($\sigma=1$).

\subsection{Problem-specific adaptation setup}

The fine-tuning and re-training datasets are generated from G-ROM solutions at $T=0.5\pi$. One-dimensional profiles are extracted, classified as DIOs or SOSs, and used to construct Gegenbauer reconstruction-error matrices over the prescribed $(m,\lambda)$ grid. To reduce redundancy, only half of the eligible profiles from each realization are retained. 
For DIOs, fine-tuning and re-training use $230$ and $2152$ samples from $6$ and $60$ initial conditions, respectively; for SOSs, the corresponding datasets contain $232$ and $2320$ samples from $4$ and $40$ initial conditions.

The training procedure follows the previous experiments, with fine-tuning performed for $200$ epochs using an initial learning rate of $8\times10^{-4}$ and re-training performed for $600$ epochs.

\section{Additional Details for 2D Inviscid Burgers' Problem\label{appendix_eg4}}

\subsection{G-ROM baseline construction\label{appendix_eg4_grom}}

The spatial domain is discretized using a $256\times256$ uniform grid. 
The reduced-order space is constructed from snapshots generated by evolving $20$ randomly sampled initial conditions from $t=0$ to $t=2$ using the method of characteristics, with $100$ time instances collected for each trajectory. 
The snapshots are then smoothed by a Gaussian filter with width $\sigma=1$ and compressed by truncated SVD, yielding $r=40$ POD modes to span the G-ROM space.
The G-ROM system evolves using Galerkin projection with third-order Runge-Kutta time stepping during the online stage.

\subsection{Problem-specific adaptation setup\label{appendix_eg4_adaptation_setup}}

For the problem-specific adaptation datasets, the G-ROM solution is evolved to $T=1.0$, and one-dimensional profiles are extracted by randomly selecting $10$ lines in the $y$-direction from each two-dimensional solution. 
This yields $150$ fine-tuning samples from $15$ initial conditions and $2204$ re-training samples from $250$ initial conditions. 
The training procedure follows the previous examples, except that fine-tuning uses a reduced initial learning rate of $5\times10^{-4}$.

\section{Additional Details for Accuracy-Efficiency Trade-off\label{appendix_trade_off}}

\subsection{Accuracy and efficiency versus fine-tuning dataset volume \label{appendix_finetuning_data_volume}}

\begin{figure}[!ht]
    \centering
    \includegraphics[width=0.75\linewidth]{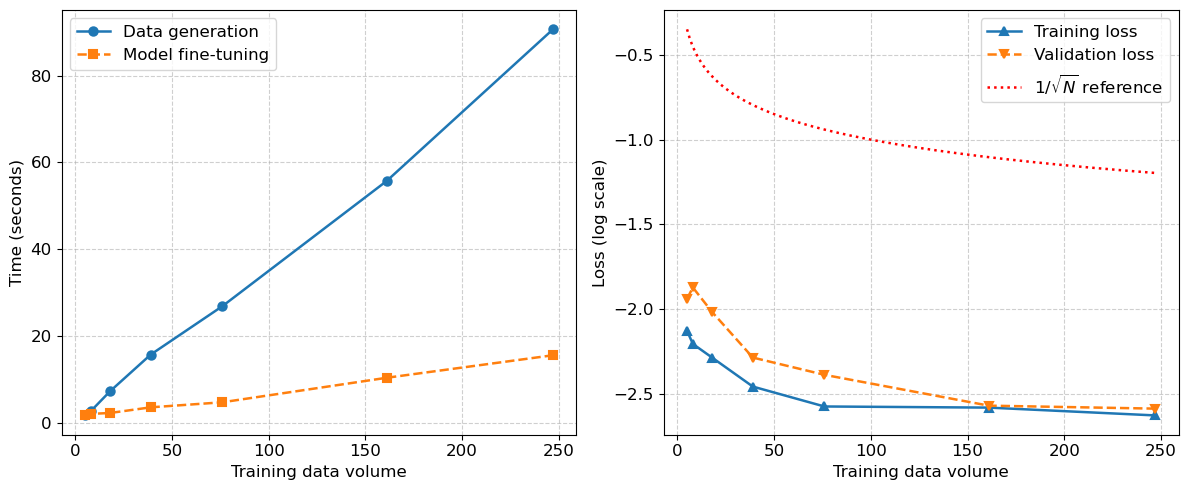}
    \caption{Impact of fine-tuning dataset volume on computational cost and training loss. Left: dataset generation and fine-tuning time as a function of dataset size. Right: optimal training and validation loss obtained from an 80/20 training/validation split.}
    \label{time_and_loss_vs_volume}
\end{figure}

Figure~\ref{time_and_loss_vs_volume} reports the computational cost and training behavior for different fine-tuning dataset sizes. 
The dataset-generation cost increases with the number of target-specific samples and dominates the fine-tuning time. 
The training and validation losses show a larger gap for small datasets, indicating overfitting, while this gap decreases as the dataset size increases.

\begin{figure}[!ht]
    \centering
    \includegraphics[width=0.98\linewidth]{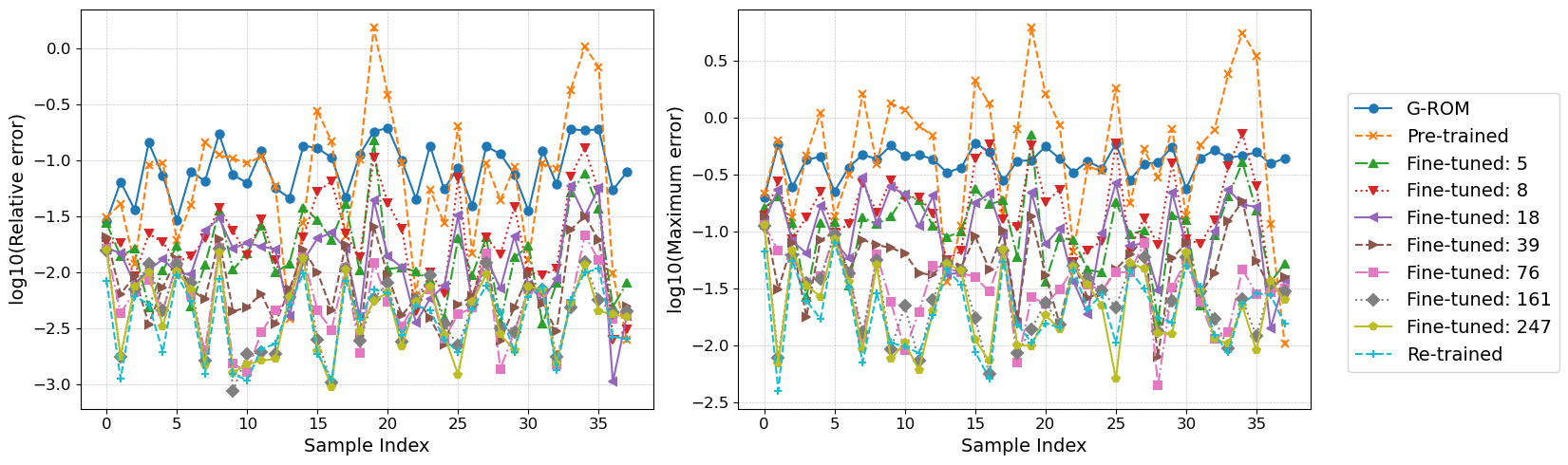}
    \caption{Reconstruction accuracy versus fine-tuning dataset volume for $38$ Burgers' G-ROM test samples. Left: relative $\ell^2$ error; right: maximum error. Comparisons include the pre-trained model, fine-tuned models with varying dataset sizes, and the fully re-trained model.}
    \label{accuracy_and_cost_vs_volume}
\end{figure}

Figure~\ref{accuracy_and_cost_vs_volume} shows the per-sample reconstruction errors over $38$ Burgers' G-ROM test samples for different fine-tuning dataset sizes, complementing the mean and quantile statistics reported in the main text.

\subsection{Accuracy and efficiency versus fine-tuning strategies}

Table~\ref{tab_finetune_grad_formula} summarizes the estimated backward-pass gradient computations required by different choices of trainable modules, explaining the empirical training-time differences observed in Table~\ref{tab_network_expressivity} and supporting the lightweight fine-tuning strategy adopted in the main text.
Here, $d_\ell$, $w_\ell$, and $k_\ell$ denote the input dimension, output width, and kernel size of a convolutional layer, respectively; $N_p$ denotes the number of spatial points; and $d_{\mathrm{fc}}$ and $w_{\mathrm{fc}}$ denote the input dimension and width of the fully-connected prediction head. The subscripts $s$, $m$, $l$, and $f$ correspond to the small-, medium-, large-receptive-field convolutional branches and the fusion convolution, respectively.

\begin{table}[!ht]
\centering
\small
\caption{Backward pass gradient computation formulas for different fine-tuning strategies.}
\label{tab_finetune_grad_formula}
\begin{tabular}{c|l|l}
\hline
\textbf{Strategy} & \textbf{Trainable Modules} & \textbf{Computed Layers (Backward)} \\
\hline
\textbf{only fc} & \texttt{fc} & $2 d_\text{fc} w_\text{fc} + 2 d_\text{fc} w_\text{fc} + w_\text{fc}$ \\
\hline
\textbf{fc+fusion} & \texttt{fc, fusion\_conv} & \makecell[l]{fusion: $2 d_f w_f k_f N_p + 2 d_f w_f k_f + w_f N_p$ \\ 
fc: $2 d_\text{fc} w_\text{fc} + 2 d_\text{fc} w_\text{fc} + w_\text{fc}$} \\
\hline
\textbf{fc+fusion+large} & \makecell[l]{\texttt{fc, fusion\_conv,} \\ \texttt{conv\_large}} & \makecell[l]{large: $2 \times (2 d_l w_l k_l N_p + 2 d_l w_l k_l + w_l N_p)$ \\
fusion: $2 d_f w_f k_f N_p + 2 d_f w_f k_f + w_f N_p$ \\
fc: $2 d_\text{fc} w_\text{fc} + 2 d_\text{fc} w_\text{fc} + w_\text{fc}$} \\
\hline
\textbf{only fusion} & \texttt{fusion\_conv} & $2 d_f w_f k_f N_p + 2 d_f w_f k_f + w_f N_p$ \\
\hline
\textbf{only conv} & \makecell[l]{\texttt{conv\_small,}\\ \texttt{conv\_medium, } \\ \texttt{conv\_large}}& \makecell[l]{small: $2 \times (2 d_s w_s k_s N_p + 2 d_s w_s k_s + w_s N_p)$ \\
medium: $2 \times (2 d_m w_m k_m N_p + 2 d_m w_m k_m + w_m N_p)$ \\
large: $2 \times (2 d_l w_l k_l N_p + 2 d_l w_l k_l + w_l N_p)$} \\
\hline
\textbf{all} & All Layers & \makecell[l]{small: $2 \times (2 d_s w_s k_s N_p + 2 d_s w_s k_s + w_s N_p)$ \\
medium: $2 \times (2 d_m w_m k_m N_p + 2 d_m w_m k_m + w_m N_p)$ \\
large: $2 \times (2 d_l w_l k_l N_p + 2 d_l w_l k_l + w_l N_p)$ \\
fusion: $2 d_f w_f k_f N_p + 2 d_f w_f k_f + w_f N_p$ \\
fc: $2 d_\text{fc} w_\text{fc} + 2 d_\text{fc} w_\text{fc} + w_\text{fc}$} \\
\hline
\end{tabular}
\end{table}

\subsection{Out-of-distribution generalization\label{appedix_ood}}

\begin{figure}[!ht]
    \centering
    \includegraphics[width=0.8\linewidth]{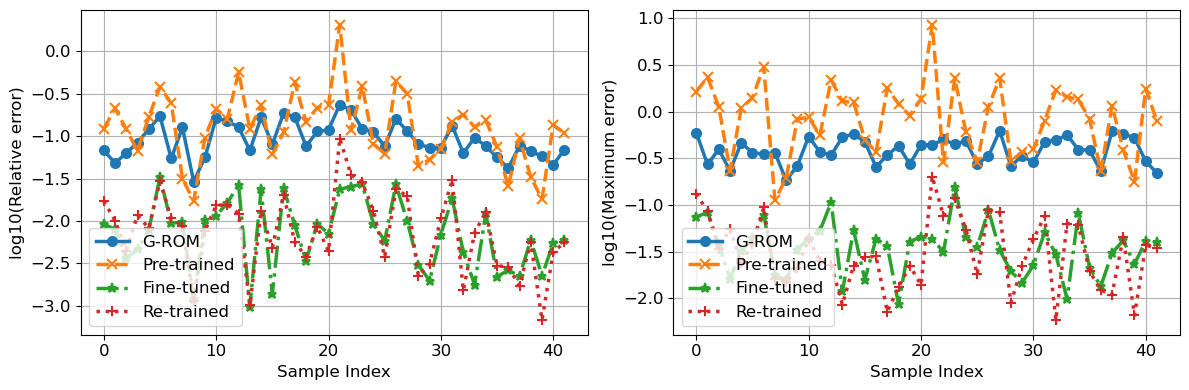}
    \caption{Out-of-distribution test on Burgers' problem with time shift: reconstruction errors over $42$ G-ROM samples. 
    Left: relative $\ell^2$ error; right: maximum pointwise error. Results are reported for pre-trained, fine-tuned, and fully re-trained models.}
    \label{time_OoD}
\end{figure}

Figure~\ref{time_OoD} reports the temporal-shift test, where Burgers' G-ROM samples are evaluated at times outside the fine-tuning target distribution.

\begin{figure}[!ht]
    \centering
    \includegraphics[width=0.8\linewidth]{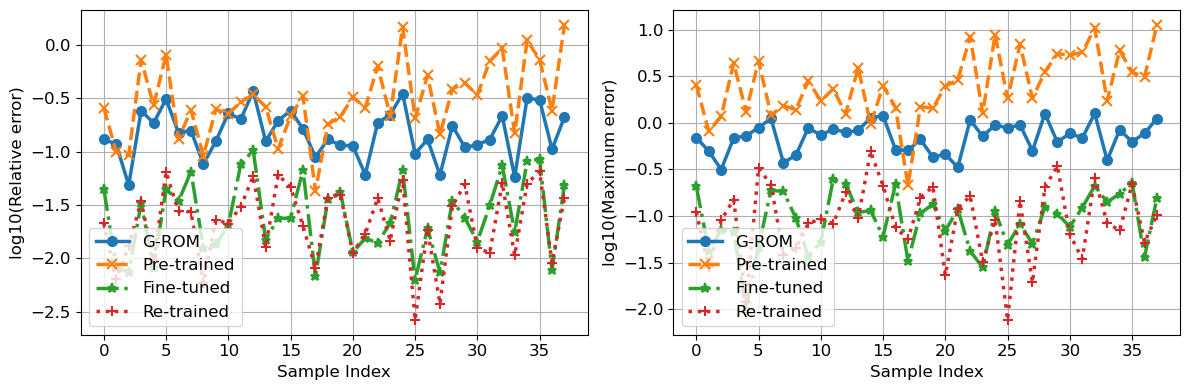}
    \caption{Out-of-distribution test on Burgers' problem with added Gaussian noise: reconstruction errors over $38$ G-ROM samples. 
    Left: relative $\ell^2$ error; right: maximum pointwise error. Results are reported for pre-trained, fine-tuned, and fully re-trained models.}
    \label{noise_OoD}
\end{figure}

Figure~\ref{noise_OoD} reports the additive-noise test, where Gaussian perturbations are added to the G-ROM inputs. 

\section*{Ethics declarations}
\noindent\textbf{Conflict of interest} \\
The authors declare that they have no conflict of interest.

\section*{Data Availability}
\noindent Data sets generated during the current study are available from the corresponding author on reasonable request.

    \bibliographystyle{abbrv}
    \bibliography{ref}

\end{document}